\documentclass[10pt]{report}
\usepackage{oldgerm}
\usepackage{amsmath,amssymb,amsfonts}
\author{}  
\usepackage{epsfig}
\usepackage{pstricks}
\usepackage{pst-node,pst-tree}
\usepackage{graphicx}
\usepackage{epic}
\usepackage{makeidx}
\usepackage{eepic}
\usepackage{delarray}
\usepackage{letterspace}
\usepackage{boxedminipage}
\newcommand{\E}{{\mathbb E}}
\newcommand{\N}{{\mathbb N}} 
\newcommand{\R}{{\mathbb R}}
\newcommand{\F}{\mathcal{F}}
\newcommand{\G}{\mathcal{G}}
\newcommand{\C}{\mathcal{C}}
\newcommand{\U}{\mathcal{U}}
\newcommand{\J}{\mathcal{J}}
\date{}
\begin{document}
\begin{center}CONVEX POLYGONS\\ AND THE ISOPERIMETRIC PROBLEM\\ IN SIMPLY CONNECTED SPACE FORMS ${\Large M_\kappa^2}$ \\[2mm]
\end{center}
\begin{center}
Anisa M. H. Chorwadwala 
and  A. R. Aithal\\[2mm]
\end{center} 
ABSTRACT.\quad In this article, we prove that {\it there exists a unique perimeter minimizer among all piecewise smooth simple closed curves in} $M^2_\kappa$ {\it enclosing area } $A>0$ ($ A \leq 2 \pi $ if $\kappa=1$), {\it and it is a circle in} $M^2_\kappa$ {\it of radius} $AS_\kappa \left(\sqrt{A \, (4 \pi - \kappa\,A)}/(2 \pi)\right)$, where $AS_\kappa(t):= \begin{cases} t & \text{if $\kappa=0$},\\
\arcsin (t) & \text{if $\kappa=1$},\\
\sinh^{-1}(t) & \text{if $\kappa =-1$}.
\end{cases} $ \\[1mm]
We also prove the isoperimetric inequality for $M^2_\kappa$. We give an elementary geometric proof which is uniform for all three simply connected space forms. \\[3mm]
\begin{center}
0. INTRODUCTION
\end{center}
\quad  Questions of the following type arise quite naturally. Why are small water droplets and bubbles that float in air approximately spherical? Why does a herd of reindeer form a circle if attacked by wolves? Of all geometric figures having certain property, which has greatest area or volume; and of all figures having certain property, which has least perimeter or surface area? These problems are capable of stimulating mathematical thought. 

The isoperimetric problem on a surface is to enclose a given area with the shortest possible curve. The classical isoperimetric theorem asserts that in the Euclidean plane the unique solution is a circle. This property of the circle is most succinctly expressed in the form of an inequality called the isoperimetric inequality. The solution of isoperimetric problem for `rectangles' was already known to Euclid. Little progress was made from Greek geometers until Swiss mathematicians Simon L'Huilier and Jacob Steiner of late eighteenth century. Using a symmetry argument Steiner has shown that the minimizer is a circle. However he did not prove the existence of a minimizer. By the use of `approximating polygons', Edler filled this gap in 1882. However, these methods have long been forgotten and seem to have been rediscovered in \cite{How-Hutch-Morg}. Here, by analogous methods, we solve the isoperimetric problem on the simply connected surface $M_\kappa^2$ having constant sectional curvature $\kappa$ ($\kappa=0, \pm 1$)
, and prove that `circle' is the unique solution to the isoperimetric problem. In this article, we give an elementary geometric proof which is uniform for all three simply connected space forms. 

Before starting, a little more history is worth inserting. The history included here is taken mainly from the survey article of Osserman \cite{Oss1} which is about developments in the theory of isoperimetric inequalities. This survey recounts some of the most interesting of the many sharpened forms, various geometric versions, generalizations, and applications of this inequality. Also see the book by H. Hadwiger \cite{Hadw1}, Other general references given in \cite{Oss1} are Kazarinoff \cite{Kazar1}, P$\acute{o}$lya [\cite{Polya2}, Chapter X], Porter \cite{Porter1}, and the books of Blaschke listed in the bibliography. One aspect of the subject is given by Burago \cite{Burago-Zalgaller1}. Others may be found in \cite{Osser4} and in the book of Santal$\grave{o}$ \cite{Santalo4}.

Most histories of the isoperimetric problem begin with its legendary origins in the ``Problem of Queen Dido''. Her problem (or at least one of them) was to enclose an optimal portion of land using a leather thong fashioned from oxhide. If Dido's was the true original isoperimetric problem, then what is needed is a solution not in the plane, but on a curved surface. (For more history of the classical case of curves  in the plane see Mitrinovi$\acute{c}$ \cite{Mitrinovic1}). The consideration of the isoperimetric problem on curved surfaces goes quite a way back, at least to an 1842 paper of Steiner \cite{Steiner1}. The fact that the smooth closed curve solving the isoperimetric problem on a surface must have constant geodesic curvature was mentioned in Steiner's paper [\cite{Steiner1}, p.150], and a proof was given in 1878 by Minding \cite{Minding1}. A detailed discussion is given in $\S 18$ of an extraordinary paper of Erhard Schmidt \cite{Schmidt4}. This paper provides an extended analysis of the isoperimetric problem on surfaces. 
 An interesting solution to the isoperimetric problem for curves on the sphere was given by F. Bernstein in 1905 \cite{Bernstein1}. A proof of the isoperimetric inequality for the hyperbolic plane was given in 1940 by Schmidt [\cite{Schmidt2}, p.209]. Fiala \cite{Fiala1} appears to have been the first to prove a general isoperimetric inequality for surfaces of variable Gauss curvature. See also Bol \cite{Bol1},  Schmidt [\cite{Schmidt4}, p.618],  Aleksandrov \cite{Alek1}, [\cite{Alek2}, p.509] and  Aleksandrov and Strel'cov \cite{Alek3, Alek4}. For the survey of the isoperimetric problem on general Riemannian manifolds refer to [\cite{Oss1}, p.1211, \S C].

 The fact is, the isoperimetric inequality holds in the greatest generality imaginable, but one needs suitable definitions even to state it. The isoperimetric inequalities have proved useful in a number of problems in geometry, analysis, and physics.

   We remark here that there are many other results of a similar nature, referred to as isoperimetric inequalities of mathematical physics, where extrema are sought for various quantities of physical significance such as the energy functional or the eigenvalues of a differential equation. They are shown to be extremal for a circular or spherical domain. Faber-Krahn Theorem \cite{Faber1, Krahn1, Krahn2} is an example of such results. Please see Rayleigh's fundamental treatise {\sf The theory of sound} [\cite{Rayleigh1}, $\S 210$]. Extensive discussions of such problems can be found in the book of P$\acute{o}$lya and Szeg$\ddot{o}$ \cite{Polya-Szego1} and the review article by Payne \cite{Payne1}. For some recent results of this type see \cite{Ramm-Shivakumar, Kesavan, Anisa-Aithal, Anisa-Vemuri, Anisa-Rajesh-Kesavan}. For specific relations between the first non-zero eigenvalue of the Laplacian and geometric isoperimetric constants associated with compact Riemannian manifold, we refer to papers of Cheeger \cite{Cheeg1} and Yau \cite{Yau1}. (See also Buser \cite{Bus1, Bus3, Bus5, Bus6}, Berger \cite{Berger1}, Chavel \cite{Chavel1} and  Reilly \cite{Reilly3}). \\[3mm]

We now state the main results :\\
\noindent \textbf{Theorem 1}\quad{\it Fix} $n \geq 3$ {\it in} $\N ~ \& \; A \in  \begin{cases} 
(0, \infty) & \text{if $\kappa=0$},\\
(0,2\pi) & \text{if $\kappa=1$},\\       
(0,(n-2)\,\pi) & \text{if $\kappa=-1$} .
\end{cases}$\\
 {\it Among all polygons with} $n$ {\it sides in} $M^2_\kappa$ {\it having area
} $A$, {\it perimeter minimizer is the regular} $n${\it -gon
}.\\[2mm]
Let $AS_\kappa(t):= 
\begin{cases} t & \text{if $\kappa=0$},\\
\arcsin (t) & \text{if $\kappa=1$},\\
\sinh^{-1}(t) & \text{if $\kappa =-1$}.
\end{cases}$\\[1.5mm] 
\noindent \textbf{Theorem 2}\quad{\it Fix } $A >0$ ($A \leq 2 \pi$ if $\kappa=1$). {\it There exists a unique perimeter minimizer among all piecewise smooth simple closed curves in} $M^2_\kappa$ {\it enclosing area} $A$, {\it and it is a circle in} $M^2_\kappa$ {\it of radius} $AS_\kappa \left(\sqrt{A \, (4 \pi - \kappa\, A)}/(2 \pi)\right)$.\\[1.5mm]
\noindent \textbf{Corollary 3}\quad {\it Fix} $A >0$ ($A \leq 2 \pi$ if $\kappa=1$). {\it There exists a unique perimeter minimizer among all piecewise smooth simple closed curves in} $M^2_\kappa$ {\it having} $m$ {\it components each enclosing area} $A_i>0$ {\it such that} $A=\sum_{i=1}^m A_i$, {\it and it is a circle in} $M^2_\kappa$ {\it of radius} $AS_\kappa \left(\sqrt{A \, (4 \pi - \kappa\, A)}/(2 \pi)\right)$.\\[1.5mm]
\noindent \textbf{Theorem 4 (The Isoperimetric Inequality for $\boldsymbol{M^2_\kappa}$)}\quad {\it For any piecewise smooth simple closed
curve } $\C$ {\it in  }
 $M^2_\kappa$ {\it with arc-length} $\ell$ {\it and enclosing area} $A>0$ ($A \leq  2 \pi$ if $\kappa=1$) 
{\it we have} $ \ell
 ^{\,2} \geq 4 \pi A - \kappa \,A^2$ {\it and equality holds if and only if } $\C$ {\it is a circle in} $M^2_\kappa$ {\it of radius} $AS_\kappa \left(\sqrt{A \, (4 \pi -\kappa \, A)}/(2 \pi)\right)$.
\par In section 1, we introduce the model spaces $M_\kappa^2$ (as Riemannian manifolds) and discuss isometries of $M_\kappa^2$. 
In sections 2, we state few results on triangles and polygons in $M^2_\kappa$ and we have given proofs mostly when the results are not available in books. 
Regular polygons in $M^2_\kappa$ are studied in section 3. Section 4 contains the proof of Theorem 1. In section 5,
proofs of Theorem 2, Corollary 3, Theorem 4 are given. Section 6 is an appendix to this article. 
\par 

\begin{center}
1. ISOMETRIES OF $M_\kappa^2$ 
\end{center}
A space form is a complete Riemannian manifold with constant sectional curvature $\kappa$. Complete, simply connected Riemannian manifolds of dimension $d$, with constant sectional curvature $\kappa $ are denoted by $M_\kappa^d$.

Let $ <~,~ >_0$ denote the standard inner product of the Euclidean space $\E^d$ ($d \in \N$). The Euclidean space $(\E^2, <~,~ >_0)$ and $S^2=\{x \in \E^3\;| <x,x>_0\,=1\}$, the unit sphere in $\E^3$ with induced Riemannian metric from $\E^3$ are the model spaces for $M_0^2$ and $M_1^2$ respectively. The hyperboloid of one sheet $ \{(x_1, x_2,x_3) \in \R^{3}\, | x_1^2+ x_2^2 -x_{3}^2= -1\;\& \;x_{3}>0\}$ with the Riemannian metric induced from the quadratic form $\left<x, y\right>_{-1}:=x_1\, y_1+x_{2}\,y_{2}-x_{3}\,y_{3}$ where $x=(x_1, x_2,x_3),\,  y=(y_1, y_2,y_3)$ is the model space for $M_{-1}^2$. The {\it inner metric} $d_\kappa$ of $M^2_\kappa$ is given by the formula 
$$d_\kappa(x, y)= 
\begin{cases}
\;\sqrt{\left<x-y, x-y\right>_0}& \text{if $\kappa=0$}\\ 
\; AC_\kappa \left(\kappa\,<x, y>_\kappa \right) & \text{if $\kappa \neq 0$}
\end{cases} ~~\;\forall \;x, y \in M_\kappa^2, $$ 
where 
$$ AC_\kappa (t):= 
\begin{cases} t & \text{if $\kappa=0$},\\
\arccos (t) & \text{if $\kappa=1$},\\
\mbox{arccosh}(t) & \text{if $\kappa =-1$}.
\end{cases}$$ 
For $p \in M^2_\kappa$ and $r > 0$ ($r<\pi$ if $\kappa=1$), \,$B_\kappa(p,r):= \{x \in M^2_\kappa \, |\, d_\kappa(p,x)<r\}$ denotes the open ball in $M^2_\kappa$ with center $p$ and radius $r$. Its boundary is $\C_\kappa(p,r):=\{x \in M_\kappa^2 \, | \,d_\kappa(p,x)=r\}$. \\[3mm]
If we take 
$$p_0= \begin{cases}
(0,0) & \text{if $\kappa=0$}\\
(0,0,1) & \text{if $\kappa \neq 0$}
\end{cases}\eqno{(1)}$$ 
then $\C_\kappa(p_0,r)$ is nothing but 
a Euclidean circle in the plane $\{\left(x_1,x_2, |\kappa|\, C_\kappa (r)\right)\, |\; x_1, x_2 \in \R\}\subseteq \R^3$ with center $ C_{\kappa} (r) \, p_0$ and radius $S_{\kappa} (r)$, where\\
 $$C_\kappa (t)= \begin{cases}
t & \text{if $\kappa=0$}\\
\cos t & \text{if $\kappa=1$}\\
\cosh t &\text{if $\kappa=-1$}
\end{cases}~~~\mbox{ and } ~~~S_\kappa (t)= \begin{cases}
t & \text{if $\kappa=0$}\\
\sin t & \text{if $\kappa=1$}\\
\sinh t &\text{if $\kappa=-1$.}
\end{cases}$$   
We say that $\C_\kappa(p,r)$ is a {\it circle} in $M_\kappa^2$ of radius $r$. The area of the ball $B_\kappa(p,r)$ is $4\, \pi\, S_\kappa^2\left( \dfrac{r}{2}\right) $. The perimeter of the ball $B_\kappa(p,r)$ is $2\, \pi\, S_\kappa\left( r \right) $.\\

Let $\tilde{H}_0$ denote a line in $\E^2$. Let $\tilde{H}_\kappa$, $\kappa \neq 0$, denote a 2-dimensional vector subspace of $\R^3$. Let $n_\kappa$ be a unit vector normal to $\tilde{H}_\kappa$ at any point of $\tilde{H}_\kappa$. 
Let $H_\kappa:= \tilde{H}_\kappa \cap M_\kappa^2$. We call $H_\kappa$ a {\it line} in $M_\kappa^2$. Then $M_\kappa^2 \setminus H_\kappa$ has two connected components. We call these components having $H_\kappa$ as common boundary as {\it open half-spaces} in $M_\kappa^2$. When $\kappa =1$ they are the open hemispheres in $S^2$.\\[2mm]
\noindent \textbf{Definition :}\quad The reflection $r_{H_\kappa}$ through a {\it line} $H_\kappa$ in $M_\kappa^2$ is defined as $$r_{H_\kappa}(x) = x-2\, <x, n_\kappa>_\kappa \, n_\kappa.$$ 
\noindent \textbf{Definition :}\quad Let $(M,g)$ be a Riemannian manifold. A diffeomorphism $\varphi :M \rightarrow M$ is called an {\it isometry} of $(M,g)$ if the differential $ d\varphi $ preserves Riemannian metric, i.e., for all $ x\in M $ and for all pairs $ u,v\in T_{x}M $ we have $$ g_{x}( u,v) = g_{\varphi(x)}( d\varphi|_{x}(u),  d\varphi|_{x}(v) ) .$$
\noindent \textbf{Remark :}\quad Any isometry $ \varphi $ of ($ M,g $) satisfies $ d(\varphi(x),\varphi(y))=d(x,y) ~\forall \,x,y  \in M $, where $d$ is the inner metric of ($M,g$).\\[2mm]
\noindent \textbf{Proposition 1.1}\quad{\it Given any positive integer} $k$ {\it and two sets of} $k$ {\it points} $\{A_1, \ldots, A_k\}$ {\it and} $\{B_1, \ldots, B_k\}$ {\it in} $M_\kappa^2$ {\it such that} $d_\kappa (A_i, A_j)=d_\kappa (B_i, B_j)$ $\forall \; i,j \in \{1,\ldots, k\}$ {\it there exists an isometry of} $M_\kappa^2$ {\it mapping} $A_i$ {\it to} $B_i$ $\forall \; i \in \{1,\ldots, k\}$. {\it Moreover, one can obtain such an isometry by composing} $k$ {\it or fewer reflections through lines}. 
(cf. \cite{Bridson}) \\[3mm]
\noindent \textbf{Proposition 1.2}\quad {\it Let} $\phi$ {\it be an isometry of} $M_\kappa^2$. 
\begin{itemize}
\item[(1)] {\it If} $\phi$ {\it is not the identity, then the set of points which it fixes is contained in a line}.
\item[(2)] {\it If} $\phi$ {\it acts as the identity on some line} $H_\kappa$, {\it then} $\phi$ {\it is either the identity or the reflection} $r_{H_\kappa}$ {\it through the line} $H_\kappa$.
\item[(3)] $\phi$ {\it can be written as the composition of three or fewer reflections through lines}.
\end{itemize} (cf. \cite{Bridson})\\[3mm]

We now describe the Isometry group of the model spaces $M_\kappa^2$, denoted as Iso($M_\kappa^2$). Let $O(d)$, $d \in \N$, denote the group of orthogonal matrices, i.e., those real $d \times d$ matrices $A$ which satisfy $^tA\, A=\mbox{Id}$, where $^tA$ is the transpose of $A$ and $\mbox{Id}$ is the identity matrix. Consider the group $GL(d+1,\R)$ (thought of as matrices) with the usual linear action on $\R^{d+1}$. Let $O(d,1)$ denote the subgroup of $GL(d+1,\R)$ 
consisting of those matrices which leave invariant the bilinear form $<\cdot,\cdot>_{-1}$. A simple calculation shows that $O(d,1)$ consists of those $(d+1)\times (d+1)$ matrices $A$ such that $^t A\,J\,A=J$, where $J$ is the diagonal matrix with entries $(1,1,\ldots,1,-1)$ in the diagonal. 
Let $O(d,1)_+ \subseteq O(d,1)$ be the subgroup consisting of those matrices in $O(d,1)$ whose bottom right hand entry is positive.\\[3mm] 
\noindent \textbf{Proposition 1.3}\quad {\it }
\begin{itemize}
\item[(i)] Iso($M_0^2$)$\cong$ $\R^2 \rtimes O(2)$, the semidirect product.
\item[(ii)] Iso($M_1^2$) $\cong O(3)$.
\item[(iii)] Iso($M_{-1}^2$) $\cong O(2,1)_+$.
\end{itemize}
(cf. \cite{Bridson})
\begin{center}
2. GEODESIC SEGMENTS, TRIANGLES AND POLYGONS IN $M^2_\kappa$ 
\end{center}
\noindent \textbf{Definition :}\quad Connected subsets of {\it line} $H_\kappa$ in $M_\kappa^2$ are called {\it geodesic segments} of $M_\kappa^2$.\\

Consider $x,y \in M_\kappa^2$ such that $x \neq y$ ($x \neq \pm y$ when $\kappa=1$). Put
 $$ v := 
y - \kappa\, <y, x>_\kappa x+ (|\kappa|-1)x
.$$ 
Then $v \in T_xM_\kappa^2$. We denote 
 $$\left\{\left.C_\kappa(t)^{|\kappa|} \, x + S_\kappa(t)\,\dfrac{v}{\sqrt{<v,v>_\kappa}}\;\right|\; 0\leq t \leq d_\kappa(x,y)\right\}~\mbox{ by }~[x,y].$$ Then $[x,y]$ is a geodesic segment in $M_\kappa^2$ joining $x$ and $y$. 

For $p \in M_\kappa^2$ and unit vector $v \in T_pM_\kappa^2 \setminus \{\underline{0}\}$, let $\gamma_{p,v}$ denote the geodesic with the initial conditions $\gamma_{p,v}(0)=p$ and $\gamma^\prime_{p,v}(0)=v$. Then
$$\gamma_{p,v}(t)= C_\kappa(t)^{|\kappa|}\, p + S_\kappa(t)\,\dfrac{v}{\sqrt{<v,v>_\kappa}}~~~~~~(t \in \R). $$
 A {\it polygon} $\wp$ in $M^2_\kappa$ is a closed region whose boundary
$\partial \wp$ is a simple closed curve (i.e., it is homeomorphic to $S^1$) consisting of geodesic segments. A point $p$ of
$\partial \wp$ is called a {\it vertex} of $\wp$ if $\partial \wp$ intersected with some disc with center $p$ consists of two radial geodesic segments which are not extensions of each other. The geodesic segments constituting $\partial \wp$ are called {\it sides} of $\wp$. For a vertex $p$ of a polygon $\wp$, let $\gamma_{p,v_1}$ and $\gamma_{p,v_2}$ denote the sides of $\wp$ having common vertex $p$. If we give positive orientation to $\partial \wp$ then the {\it angle} of polygon $ \wp $ at vertex $ p $ is defined as 
\begin{eqnarray*} 
\measuredangle \, \mbox{at} \, p := \begin{cases} 
~~~\measuredangle{ \{v_1,v_2 \} }~
~~~\hfill\mbox{if det}\left({v_1,v_2}\right)< 0, \\[1mm]
 ~2\pi\,-\,\measuredangle{ \{v_1,v_2
\} }~~~~\hfill \mbox{if det}\left({v_1,v_2
}\right) > 0.
\end{cases}
\end{eqnarray*}

A polygon $\wp$ is said
to be {\it convex} if for any $x,y \in \wp$ (with $y \neq -x$ if $\kappa =1$), the geodesic segment $[x,y]$ is contained in $\wp$.
A polygon $\wp$ is said
to be {\it locally convex} if for any $x \in \wp$, $B_\kappa(x,r) \cap \wp$ is convex $\forall ~ r >0$. Note that a connected locally convex polygon is convex and vice versa. A polygon in $M_1^2$ is called {\it proper polygon}
if it contains no pair of antipodal points. A polygon (proper polygon if $\kappa=1$) of $n$ sides is called an {\it n}{\it -gon} in $M_\kappa^2$. Note that for any $n$-gon, $n \geq 3$ always holds. For $\kappa \neq 1$, any 3-gon is always convex. A convex $3$-gon in $ M_\kappa^2 $ is called a {\it triangle} in $ M^2_\kappa $. A triangle in $ M_\kappa ^2 $ having vertices $ x,\, y,\, z \in M_\kappa ^2 $ is denoted by $ [x,y,z] $. \\[2mm] 
\noindent {\it \textbf{Law of Cosine}} for triangles in $M_\kappa^2$ :\\[1mm]
\noindent \underline{$\kappa=0$} \quad $c^2= a^2 + b^2 - 2\, a\, b\, \cos \gamma$, \\
\noindent \underline{$\kappa \neq 0$} \quad 
$C_\kappa (c) =  C_\kappa (a)~C_\kappa (b) + \kappa \, S_\kappa (a)~S_\kappa (b) ~\cos \gamma$,\\
where $a$, $b$, $c$ are the sides of the triangle and $\gamma$ is the angle opposite to side $c$. \\[2mm]
In particular, fixing $a$, $b$ and $\kappa$, one sees that $c$ is a strictly increasing function of $\gamma \in [0,\pi]$. The triangle inequality for a triangle in $M_\kappa ^2$ follows from the Law of Cosine. Strict triangle inequality holds for triangles in $M_\kappa ^2$.\\[2mm] 
\noindent {\it \textbf{Law of Sine}} for triangles in $M_\kappa^2$ :\\[1mm]
$$\dfrac{S_\kappa(a)}{\sin \alpha}=\dfrac{S_\kappa(b)}{\sin \beta}=\dfrac{S_\kappa(c)}{\sin \gamma},$$
where $a$, $b$, $c$ are the sides of the triangle and $\alpha$, $\beta$, $\gamma$ are the angles opposite to sides $a$, $b$, $c$ respectively. \\[3mm]
\noindent \textbf{Theorem 2.1}\quad{\it The area of a triangle} $T$ {\it in} $M^2_\kappa$ ($\kappa \neq 0$) {\it with angles} $\alpha,\, \beta, \,\gamma $ {\it is equal to}
$\kappa \, \left(\alpha + \beta + \gamma - \pi\right)$.\\[1mm]
\noindent {\it Proof.}\quad By Gauss-Bonnet Formula, $\left(\alpha + \beta + \gamma - \pi\right)$ is nothing but $\int_T \kappa \, dV$, where $dV$ is the area element of $M_\kappa^2$. Therefore, for $M_\kappa^2$ ($\kappa \neq 0$), area of triangle $T$ is equal to $\kappa \, \left(\alpha + \beta + \gamma - \pi \right)$.\hfill $\blacksquare$\\[1.5mm]
\noindent \textbf{Remarks :}\quad 
\begin{itemize}
\item[(1)] For $\kappa=1$, Theorem 2.1 is known as {\it Girard's Theorem}. 
\item[(2)] The area of a triangle in $\E^2$ can not be determined only from its three angles.
\item[(3)] The area of the disk $B:=B_{-1}\left(p, 2\, \sinh^{-1} \left( \frac{1}{2}\right)\right)$, $p \in M_{-1}^2$, is $\pi $ which is greater than area of any triangle in $ M_{-1}^{2}$. Hence there is no triangle in  $M_{-1}^{2} $ which can inscribe the disk $B$. \textbf{Triangles in} $\boldsymbol{M_{-1}^{2}} $ \textbf{are thin !}\\
\end{itemize}
\noindent \textbf{Theorem 2.2}\quad{\it  The area} $A$ {\it of a triangle in} $M^2_\kappa$ {\it with sides} $a,b,c$ {\it is given by the equation}
$$ T_{|\kappa|} (A/4) = \sqrt{T_\kappa (s/2)\; T_\kappa \left[(s-a)/2\right] \; T_\kappa \left[(s-b)/2\right] \; T_\kappa \left[(s-c)/2\right]}\eqno{(2)}$$ {\it where} $s := (a+b+c)/2\,$ {\it and} $T_\kappa(t):= \begin{cases} t & \text{if $\kappa=0$},\\
\tan t & \text{if $\kappa=1$},\\
\tanh t & \text{if $\kappa =-1$}.
\end{cases} $ 
\\
\noindent {\it Proof.}\quad \\
\noindent \textbf{\underline{$\kappa = 0$} :}\quad Let $\gamma$ be the angle included between the sides $a$ and $b$. From the Law of Cosine we have
$$\cos \gamma = \dfrac{a^2+b^2-c^2}{2\, a\, b}.$$
Hence, $\sin \gamma = \sqrt{1-\cos^2 \gamma}=\dfrac{2}{ a\, b}\,\sqrt{s \,(s-a)\,(s-b)\,(s-c)}$ (and the Law of Sine follows immediately). Therefore,
$$A= \dfrac{1}{2}\,a\,b\, \sin \gamma = \sqrt{s \,(s-a)\,(s-b)\,(s-c)}.$$ 
\noindent \textbf{\underline{$\kappa \neq 0$} :}\quad In what follows the equations (A-1), (A-2), (A-3), $\ldots$ and (B-1), (B-2), (B-3), $\ldots$ refer to equations from Appendix A and Appendix B respectively which appear at the end of the article. By Theorem 2.1,\begin{eqnarray*}
\tan \left( \dfrac{A}{4} \right) 
& = & \tan \left( \frac{\kappa\, (\alpha + \beta +\gamma - \pi)}{4}\right)
=\dfrac{\sin \left( \frac{\kappa\, (\alpha + \beta +\gamma - \pi)}{4}\right)}
{\cos \left( \frac{\kappa\,(\alpha + \beta +\gamma - \pi)}{4}\right)}\\[1mm]
& = & \kappa\,\frac{\sin \left( \frac{ \alpha + \beta +\gamma - \pi}{4}\right)}
{\cos \left( \frac{\alpha + \beta +\gamma - \pi}{4}\right)}~~~~~~~~~~~~~~~~~~~~~~~~~~~~~~~~~~~~~~~~~~~~~~~~~~~~\mbox{[by (A-1)]}\\[1mm]
& = & \kappa\,\frac{\sin \left( \frac{\alpha + \beta}{2}\right) - \sin \left( \frac{\pi - \gamma}{2}\right)}
{\cos \left( \frac{\alpha + \beta}{2}\right) + \cos \left( \frac{\pi - \gamma}{2}\right)} ~~~~~~~~~~~~~~~~~~~~~~~~~\mbox{[by (A-14) and (A-15)]}\\[1mm] 
&=&  \kappa\,\frac{\sin \left( \frac{\alpha + \beta}{2}\right) - \cos \frac{\gamma}{2}}
{\cos \left( \frac{\alpha + \beta}{2}\right) + \sin \frac{\gamma}{2}} ~~~~~~~~~~~~~~~~~~~~~~~~~~~~~~\,~~~~\mbox{[by (A-3) and (A-7)]}\\[1mm]
& = & \kappa\,\frac{\left[ \frac{C_\kappa\left(  \frac{a-b}{2}\right) }{C_\kappa\left(  \frac{c}{2}\right) }- 1\right]
\; \cos \frac{\gamma}{2}}{\left[ \frac{C_\kappa\left(  \frac{a+b}{2}\right) }{C_\kappa\left(  \frac{c}{2}\right) } + 1\right] \;
\sin \frac{\gamma}{2}}\, \qquad \qquad~~~~~~~~~~~~~~~~~~~
\mbox{[by (B-4) and (B-6)]}\\[1mm]
\end{eqnarray*}
Therefore,
\begin{eqnarray*}
\tan \left( \dfrac{A}{4} \right)
& = & \kappa\,\frac{C_\kappa\left(  \frac{a-b}{2}\right)  - C_\kappa\left(  \frac{c}{2}\right) }{C_\kappa\left(  \frac{a+b}{2}\right) + C_\kappa \left(  \frac{c}{2}\right) }
~\cdot~ \frac{ \cos \frac{\gamma}{2}}{ \sin \frac{\gamma}{2}}
 =  \frac{S_\kappa \left( \frac{s-a}{2}\right) S_\kappa \left(  \frac{s-b}{2}\right) }{C_\kappa\left(  \frac{s}{2} \right) \, C_\kappa\left(  \frac{s-c}{2}\right) }~\cdot~ \frac{ \cos \frac{\gamma}{2}}{ \sin \frac{\gamma}{2}}\\
& &~~~~~~~~~~~~~~~~~~~~~~~~~~~~~~~~~~~~~~~~~~~~~~~~~\mbox{[by (A-15), (A-16) and (A-1)]}\\
& = & \frac{S_\kappa \left( \frac{s-a}{2}\right) S_\kappa \left(  \frac{s-b}{2}\right) }{C_\kappa\left(  \frac{s}{2} \right) \, C_\kappa\left(  \frac{s-c}{2}\right) }~
\sqrt{\frac{S_\kappa(s) \, S_\kappa(s-c)}{S_\kappa (s-a) \, S_\kappa (s-b)}}~~~~~~~\mbox{[by (B-1) and (B-2)]}\\[1mm]
& = & \frac{S_\kappa \left( \frac{s-a}{2}\right) S_\kappa \left(  \frac{s-b}{2}\right) }{C_\kappa\left(  \frac{s}{2} \right) \, C_\kappa\left(  \frac{s-c}{2}\right) }~
\sqrt{\frac{S_\kappa\left(\frac{s}{2}\right)\,C_\kappa\left(\frac{s}{2}\right)  \, S_\kappa\left(\frac{s-c}{2}\right)\,C_\kappa\left(\frac{s-c}{2}\right) }{S_\kappa \left(\frac{s-a}{2}\right) \, C_\kappa\left(\frac{s-a}{2}\right)\,S_\kappa\left(\frac{s-b}{2}\right)\,C_\kappa\left(\frac{s-b}{2}\right)}}\\
& & ~~~~~~~~~~~~~~~~~~~~~~~~~~~~~~~~~~~~~~~~~~~~~~~~~~~~~~~~~~~~~~~~~~~~~~~~~~~\mbox{[by (A-4)]}\\[1mm]
& = & \sqrt{T_\kappa \left(\frac{s}{2}\right) T_\kappa \left(\frac{s-a}{2}\right) T_\kappa \left(\frac{s-b}{2}\right) T_\kappa \left(\frac{s-c}{2}\right)}.~~~~~~~~~~~~~~~~~~~~~~\,~~\blacksquare
\end{eqnarray*}
\noindent \textbf{Remark :}\quad Equation (2) is known as {\it Heron's formula} and {\it L'Huilier's formula} for $\kappa=0$ and $\kappa=1$ respectively.\\[2mm]
\noindent \textbf{Proposition 2.3}\quad {\it Given two sides } $a, b$ {\it  and the included angle}
$\gamma$ {\it of a 
triangle in} $M^2_\kappa$, {\it its area} $A$ {\it is given by the formula }$$ CT_{|\kappa|} (A/2) = \dfrac{CT_\kappa (a/2) \;CT_\kappa(b/2) \;(\sin^2 \gamma)^{1-|\kappa|}+ \kappa \, \cos \gamma}{\sin \gamma},$$ {\it where} $CT_\kappa(t):= \begin{cases}
t & \text{if $\kappa=0$},\\
\cot t & \text{if $\kappa=1$},\\
\coth t & \text{if $\kappa =-1$}.
\end{cases}$
\\[1mm]
\noindent {\it Proof.}\quad \\
\noindent \textbf{\underline{$\kappa=0$} :}\quad $A= \dfrac{1}{2}\, a\, b\, \sin \gamma$.\\
\noindent \textbf{\underline{$\kappa\neq 0$} :}\quad Let $\alpha, \beta $ be the other two angles of the triangle opposite to sides $a,b$ respectively. By Proposition 2.1,
\begin{eqnarray*}
\sin\dfrac{A}{2} 
&=& \sin \left( \dfrac{\kappa\,(\alpha+\beta+\gamma-\pi)}{2}\right)= \kappa\, \sin \left( \dfrac{\alpha+\beta+\gamma-\pi}{2}\right) ~~~~~~~~~~\,~~\mbox{[by (A-1)]}\\
&=& -\kappa\, \cos \left( \dfrac{\alpha+\beta+\gamma}{2}\right) ~~~~~~~~~~~~~~~~~~~~~~~~~~~~~~~~~~~~~~~~~~~~~~~~~~~~\mbox{[by (A-3)]}\\
&=& -\kappa\,\left[ \cos \left( \dfrac{\alpha+\beta}{2}\right)\,\cos \left( \dfrac{\gamma}{2}\right)-\sin \left( \dfrac{\alpha+\beta}{2}\right) \sin \left( \dfrac{\gamma}{2}\right)\right] ~~~~~~~~~~~\mbox{[by (A-6)]}\\
&=&-\kappa\,
\dfrac{C_\kappa\left(\dfrac{a+b}{2}\right)-
C_\kappa\left(\dfrac{a-b}{2}\right)}{C_\kappa\left( \dfrac{c}{2}\right) }\,  \sin\left( \dfrac{\gamma}{2}\right) \, \cos\left( \dfrac{\gamma}{2}\right) ~~~\mbox{[by (B-4) and (B-6)]}\\
&=& \dfrac{\sin \gamma}{C_\kappa\left( \dfrac{c}{2}\right)} S_\kappa\left(\frac{a}{2} \right) \, S_\kappa\left(\frac{b}{2} \right)~~~~~~~~~~~~~~~~~~~~~~~~~~~~~~~~~\mbox{[by (A-4) and (A-16)]}.
\end{eqnarray*}
Hence, $$\sin\left( \dfrac{A}{2}\right) = \dfrac{S_\kappa\left( \dfrac{a}{2}\right) \,S_\kappa\left( \dfrac{b}{2}\right)\, \sin \gamma}{C_\kappa\left( \dfrac{c}{2}\right) }
.\eqno(3)$$
\begin{eqnarray*} 
\cos\left( \dfrac{A}{2}\right) 
&=& \cos\left( \dfrac{\kappa\,(\alpha + \beta + \gamma-\pi)}{2}\right)  = \cos\left( \dfrac{\alpha + \beta + \gamma-\pi}{2}\right)~~~~~\,\,\,~~~~~\mbox{[by (A-1)]}\\
&=& \sin\left(\dfrac{\alpha + \beta + \gamma}{2}\right) ~~~~~~~~~~~~~~~~~~~~~~~~~~~~~~~~~~~~~~~~~~~~~~~~~~~~~~\mbox{[by (A-7)]}\\
&=& \sin\left(\dfrac{\alpha + \beta }{2}\right)  \cos\left( \dfrac{\gamma}{2}\right)  + \cos\left(\dfrac{\alpha + \beta }{2}\right) \sin\left(  \dfrac{\gamma}{2}\right)~~~~~\,~~~~~~\,\,~~~~\mbox{[by (A-2)]}\\
&=&  \dfrac{C_\kappa\left(\dfrac{a-b}{2}\right)}{C_\kappa\left( \dfrac{c}{2}\right)} \cos^2\left( \dfrac{\gamma}{2}\right)  + \dfrac{C_\kappa \left(\dfrac{a+b}{2}\right)}{C_\kappa\left( \dfrac{c}{2}\right) }\sin^2\left( \dfrac{\gamma}{2}\right) ~\mbox{[by (B-4) and (B-6)]}\\
&=& \dfrac{\cos^{2}\left( \dfrac{\gamma}{2}\right) \,\left[ {C_\kappa\left( \dfrac{a}{2}\right) \,C_\kappa\left( \dfrac{b}{2} \right) +\kappa\, S_\kappa\left( \dfrac{a}{2}\right) \,S_\kappa\left( \dfrac{b}{2}\right) }\right]}{C_\kappa\left( \dfrac{c}{2}\right) }~~~~~~~~~~\mbox{[by (A-7)]}\\
& & +~\dfrac {\sin^{2}\left( \dfrac{\gamma}{2}\right) \,\left[ C_\kappa\left( \dfrac{a}{2}\right) \,C_\kappa\left( \dfrac{b}{2}\right) - \kappa\,S_\kappa\left( \dfrac{a}{2}\right) \,S_\kappa\left( \dfrac{b}{2}\right) \right]}{C_\kappa\left( \dfrac{c}{2}\right) }~~~~~~~\mbox{[by (A-6)]}\\
&=& \dfrac{
C_\kappa\left( \dfrac{a}{2}\right) \,C_\kappa\left( \dfrac{b}{2} \right) +\kappa\, S_\kappa\left( \dfrac{a}{2}\right) \,S_\kappa\left( \dfrac{b}{2}\right) \,\left[ \cos^{2}\left( \dfrac{\gamma}{2}\right) -\sin^{2}\left( \dfrac{\gamma}{2}\right)\right]}{C_\kappa\left( \dfrac{c}{2}\right)} .
\end{eqnarray*}
Hence by (A-8),$$\cos\dfrac{A}{2} = \dfrac{C_\kappa\left( \dfrac{a}{2}\right) \,C_\kappa\left( \dfrac{b}{2} \right) +\kappa\, S_\kappa\left( \dfrac{a}{2}\right) \,S_\kappa\left( \dfrac{b}{2}\right) \,\cos \gamma}{C_\kappa\left( \dfrac{c}{2}\right)}.\eqno(4)$$
From (3) and (4) we get,
$$\cot\left( \dfrac{A}{2}\right) = \dfrac{CT_\kappa\left( \dfrac{a}{2}\right) \,CT_\kappa\left( \dfrac{b}{2}\right) +\kappa\, \cos \gamma}{\sin \gamma}. \eqno{\blacksquare} $$
\\[1mm]
\noindent \textbf{Definition :}\quad Let $T:=[P,Q,R], ~T^\prime:=[P^\prime,Q^\prime, R^\prime]$ be triangles in $M_\kappa^2$. We say that the triangle $T$ is {\it congruent} to $T^\prime$ if there exists an isometry $f$ of $M_\kappa^2$ such that $f(P)=P^\prime$, $f(Q)=Q^\prime$ and $f(R)=R^\prime$.\\[2mm]
\noindent \textbf{Proposition 2.4}\quad {\it Let} $T 
, ~T^\prime
$ {\it be triangles in} $M_\kappa^2$. {\it Let} $a,b,c$ ({\it resp.} $a^\prime, b^\prime, c^\prime$) {\it be the sides of} $T$ ({\it resp.} $T^\prime$). {\it Let} $\alpha, \beta, \gamma$ {\it be angles of} $T$ {\it opposite to sides} $a,b,c$ {\it respectively}. {\it Let} $\alpha^\prime, \beta^\prime, \gamma^\prime$ {\it be angles of} $T^\prime$ {\it opposite to sides} $a^\prime,b^\prime,c^\prime$ {\it respectively. Then, the following are equivalent} : 
\begin{itemize}
\item[(i)] $T$ {\it is congruent to} $T^\prime$.
\item[(ii)] $a=a^\prime$, $b=b^\prime$, $c=c^\prime$.
\item[(iii)] $\alpha=\alpha^\prime$, $b=b^\prime$, $c=c^\prime$.
\item[(iv)] $\alpha=\alpha^\prime$, $\beta=\beta^\prime$, $c=c^\prime$.
\end{itemize}
{\it Each of the above imply} 
\begin{itemize}
\item[(v)] $\alpha=\alpha^\prime$, $\beta=\beta^\prime$, $\gamma=\gamma^\prime$.
\end{itemize}
{\it For} $\kappa \neq 0$, {\it all the five statements above are equivalent}.\\[3mm]
\noindent {\it Proof.}\quad See Appendix C for a proof of this Proposition. \hfill $\blacksquare$\\[3mm]
\noindent \textbf{Proposition 2.5}\quad {\it Among all triangles in} $M^2_\kappa$ {\it whose two sides are of length}
$a,b$ ($a+b< \pi$ {\it if} $\kappa =1$), {\it area maximizer is the triangle whose vertices lie on a circle having
the midpoint of the `remaining side' as its center.}\\[1mm]
\noindent {\it Proof.}\quad 
Let $\U$ denote the family of all triangles in $M_\kappa^2$ whose two sides are
of length $a,b$ ($a+b< \pi$ if $\kappa =1$). Let $r_0:= a + b$.\\
\noindent {\it Existence} \quad Upto congruence all triangles in $\U$ lie inside $B_\kappa(p,r_0)$ where $p \in M_\kappa^2$. Since $\overline{B_\kappa(p,r_0)}$ is compact in $(M_\kappa^2, d_\kappa)$ 
 there exists an `area maximizer' $T_0$ in $\U$.\\
Let $\gamma := \gamma(T)$ be the angle of a triangle $T$ in $\U$ included between the sides of length $a, \, b$. Put
$A_\gamma = area(T)$. By Proposition  2.3, when $\kappa =0$
$$A_\gamma= \dfrac{1}{2}\, a \, b\, \sin \gamma \leq \dfrac{1}{2}\, a \, b\, \sin \dfrac{\pi}{2}= A_{\frac{\pi}{2}}.$$ So, $area(T)$ is maximum when $\gamma = \dfrac{\pi}{2}$.\\
Consider 
$\kappa \neq 0$. By Proposition  2.3,
$$ \cot (A_\gamma/2) = \dfrac{CT_\kappa (a/2) \;CT_\kappa(b/2) + \kappa\, \cos \gamma }{\sin \gamma}. \eqno{(5)}$$
Consider the unit circle $S^1$ in $\E^2$ with center at $(0,0) =:O$. Let $Q = Q(\gamma)$ be a point in $\E^2$
such that $\parallel Q - O \parallel_{\E^2} = CT_\kappa (a/2) \, CT_\kappa (b/2)$. The information on $a,b$ 
implies that $CT_\kappa (a/2) \, CT_\kappa (b/2) > 1$, and hence $Q$ lies `outside'
$S^1$. Extend the line segment $[Q,O]$ and intersect $S^1$ at $R$. Let $P=P(\gamma)$ be the point on $S^1$ such that
$\measuredangle POR = (1- \kappa)\, \dfrac{\pi}{2}+ \kappa\,\gamma$. Let $N$ be the orthogonal projection of $P$ on the line joining $Q \; \& \; R$.

\vspace{1.5cm}
 \begin{center}
\psset{unit=1mm}
\begin{pspicture}(-7,-5)(20,15)
\rput(13,1)
{
\includegraphics[height=5cm]{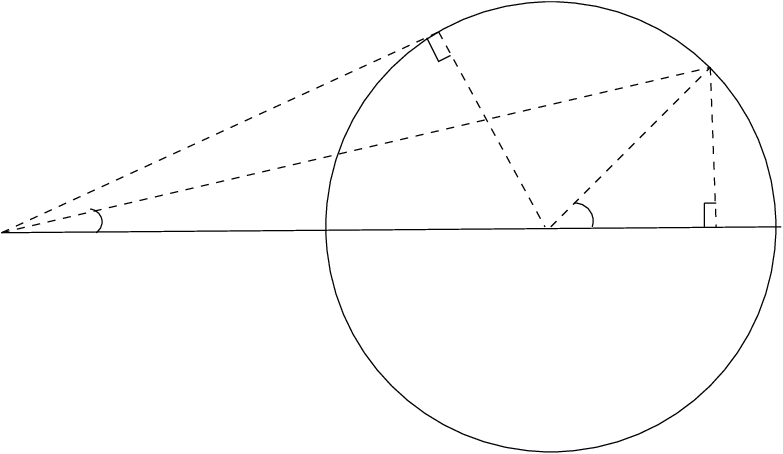}}
\rput(17.5,25){$M$}
\rput(49.5,21){$P$}
\rput(-11.5,2.5){$A_\gamma/2$}
\rput(37,3){$\gamma$}
\rput(-32,-2){$Q$}
\rput(30,-2){$O$}
\rput(48,-2){$N$}
\rput(58,0){$R$}
\rput(56,-9){\large $S^1$}
\rput(30,-30){\large Fig. 1}
\end{pspicture}
\end{center}

\vspace{3cm}
 Then,
\begin{equation*}
\parallel Q - N \parallel_{\E^2} =
\begin{cases}
\parallel Q - O \parallel_{\E^2} + \kappa \,\parallel N - O \parallel_{\E^2} & \text{if \,$\gamma \in (0, \frac{\pi}{2}]$},\\[1mm]
\parallel Q - O \parallel_{\E^2} - \kappa \, \parallel N - O \parallel_{\E^2} & \text{if \,$\gamma \in [\frac{\pi}{2},\pi)$}.
\end{cases}
\end{equation*}
Therefore, \begin{eqnarray*}
\parallel Q - N \parallel_{\E^2} &=& CT_\kappa (a/2) \, CT_\kappa (b/2) + \cos \left( (1- \kappa)\, \dfrac{\pi}{2}+ \kappa\,\gamma \right)\\
&=& CT_\kappa (a/2) \, CT_\kappa (b/2) + \kappa \,\cos \gamma .
\end{eqnarray*} 
By $(5)$, $$\measuredangle PQR = A/2 \;. \eqno{(6)}$$
\par Let $M =M(\gamma)$ be the point on $S^1$ such that line $QM$ is tangent to $S^1$ at $M$ and $M$
lies on the same side of line $QR$ as $P$. Then $\measuredangle PQR$  is maximum when $P= M$. Hence,
it follows by $(6)$ that $P(\gamma_0) = M(\gamma_0)$ where $\gamma_0 := \gamma(T_0)$. Thus, $\measuredangle POR$
being an external angle of a triangle in $\E^2$,
$$ (1- \kappa)\, \dfrac{\pi}{2}+ \kappa\,\gamma_0=\measuredangle P(\gamma_0)OR = \measuredangle M(\gamma_0)OR= \frac{\pi}{2} + \measuredangle M(\gamma_0)QO = \frac{\pi}{2} + \frac{A_0}{2} \eqno{(7)}$$
where $A_0 := area(T_0)$.
\par Let $\alpha_0, \, \beta_0$ be the angles of $T_0$ other than $\gamma_0$. As $A_0= \kappa \, \left(\alpha_0 + \beta_0 + \gamma_0 - \pi\right) $,
 $(7)$ implies that $\gamma_0 = \alpha_0 + \beta_0 $. Let $A, \, B, \, C$ be the vertices of $T_0$ having angles $\alpha_0, \, \beta_0, \, \gamma_0$ respectively. As $\gamma_0 > \alpha_0$, there is a unique point $D$ on
the side $[A,B]$ of $T_0$ such that $\measuredangle ACD = \alpha_0$.
Then $\measuredangle BCD = \beta_0$. Thus the triangles $[A,D,C] ~ \& ~[B,D,C]$ are isosceles triangles.
Hence the geodesic segments $ [A,D], \, [D,C], \, [D,B]$ are all of same length. Thus the vertices of $T_0$ lie on a circle
whose center is the midpoint of side $[A,B]$. \hfill $\blacksquare$\\[2mm]
\noindent \textbf{Proposition 2.6}\quad{\it Given} $a>0$ ($a< \pi$ if $\kappa=1$) {\it and} $\alpha \in (0, \pi)$ {\it there exists an isosceles triangle in} $M_\kappa^2$ {\it with base} $a$ {\it and base angles} $\alpha$ {\it if and only if} $\alpha \in (0, \alpha_{\kappa,a})$, {\it where} $\alpha_{\kappa,a}:= 
\begin{cases} 
\pi/2 & \text{if $\kappa=0$,}\\
\pi & \text{if $\kappa=1$,}\\ 
\arccos \left( \tanh \left( \frac{a}{2}\right) \right)& \text{if $\kappa=-1$}.
\end{cases} $\\[1mm]
\noindent {\it Proof.}\quad 
We give the proof for $\kappa = -1$. The proof for $\kappa \in \{0,1\}$ is similar and simpler.\\
\indent From Theorem 2.1 it follows that an isosceles triangle in $M_{-1}^2$ with base angles $\alpha$ exists only if $\alpha \in (0, \pi/2)$. Therefore, we consider  
$$a>0~~~\mbox{ and }~~~\alpha \in (0, \pi/2).\eqno{(8)}$$ 
Let $p_0= (0,0,1) \in M_{-1}^2$. The geodesic segment $$\gamma(t)=\gamma_{p_0,e_1}(t)= (\sinh(t), 0, \cosh(t)), ~t \in [0,a],$$ joins $p_0$ to $q:= (\sinh(a), 0, \cosh(a))$ and has length $a$. The vector 
$$v_1:=(\cos \alpha, \sin \alpha, 0)\in T_{p_0}\left( M_{-1}^2 \right)$$ 
makes an angle $\alpha$ with $e_1$ in $T_{p_0}\left( M_{-1}^2 \right)$. Let 
$$n:= \gamma^\prime\left( \frac{a}{2}\right) ,~~~ \tilde{H}:= \{x \in \R^3\;|\; <x,n>_{-1}=0\}~~~\mbox{ and } ~~~ H:=  \tilde{H}\cap M^2_{-1}.$$ 
Let $r_H$ denote the reflection in $M_{-1}^2$ through $H$. Let $$v_2:= {d(r_H)}_{p_0}(v_1)=r_H(v_1)= (-\cosh a\, \cos \alpha, \sin \alpha, -\sinh a\, \cos \alpha ).$$ Clearly, $v_2$ makes an angle $\alpha$ with $- \gamma^\prime(a)$ in $T_q\left( M_{-1}^2 \right)$. Consider the geodesics $\gamma_1= \gamma_{p_0, v_1}$ and $\gamma_2= \gamma_{q, v_2}$ of $M_{-1}^2$. Then, $$\gamma_1(t)= (\sinh t \, \cos \alpha, \sinh t \, \sin \alpha, \cosh t )$$ and  
$$\gamma_2(t)= (\cosh t \, \sinh a -\sinh t \, \cosh a \, \cos \alpha, \sinh t \, \sin \alpha, \cosh t \, \cosh a -\sinh t \, \sinh a \, \cos \alpha).$$
So $\gamma_1(t)= \gamma_2(t)$ for some $t \in \R \setminus \{0\}$ if and only if 
\begin{equation*} 
\left.
\begin{aligned}
 &\sinh t \, \cos \alpha =\cosh t \, \sinh a -\sinh t \, \cosh a \, \cos \alpha \\
&~~~~~~~~~~~~~~~~~~~~~~~\mbox{and}\\
&~~~~~\cosh t= \cosh t \, \cosh a -\sinh t \, \sinh a \, \cos \alpha. 
\end{aligned}
~~~ \right\} 
\end{equation*}
That is,
\begin{equation*}
\left.
\begin{aligned}
&~~\sinh t \, \cos \alpha\,(1+\cosh a) =\cosh t \, \sinh a \\
&~~~~~~~~~~~~~~~~~~~~\mbox{and}\\
&-\cosh t (1-\cosh a)= \sinh t \, \sinh a \, \cos \alpha. 
\end{aligned}
~~~ \right\} \eqno{(9)}
\end{equation*}
Using (A-9), (A-10) and (A-4) it is easy to see that each of the equations in (9) is equivalent to
$$\cos \alpha = \coth t \; \tanh \left(\frac{a}{2}\right).\eqno(10)$$
From (8) and (10) we get $t>0$. Now, for all $t>0$, 
$$\tanh (a/2) \, \coth t \in ( \tanh (a/2),\coth t ) \subset (0, \infty)$$ since $\tanh (a/2) \in (0,1)$ and $\coth t \in (1,\infty)$ $\forall \; t>0$. Therefore, from (10) it follows that $\cos \alpha \in (0, 1) \cap  (\tanh (a/2) , \coth t)= (\tanh (a/2), 1)$. Thus an isosceles triangle with base $a$ and base angles $\alpha$ exists if and only if $\alpha < \arccos \left(\tanh (a/2)\right)$. \hfill $ \blacksquare $\\[3mm]
\noindent \textbf{Proposition 2.7}\quad{\it Given} $0<a<s$ ($< \pi$ {\it if} $\kappa =1$) {\it there exists an isosceles triangle in} $M_\kappa^2$ {\it with base} $a$ {\it and equal sides} $s-\dfrac{a}{2}$.\\[1mm]
\noindent {\it Proof.}\quad Let $f(\kappa,a,s):= \dfrac{T_\kappa(a/2)}{T_\kappa\left(s-\frac{a}{2}\right)}$. Then, 
\begin{equation*}
\begin{aligned}
& f(0,a,s)\in (0,1) &\mbox{ since }~s-a/2 > a/2,~~~~~~\\
& f(1,a,s)\in (-1,1) &\mbox{ since }~0<a<s<\pi~~~~~~~\\
& \mbox{and}\\
& f(-1,a,s)> \tanh(a/2)\in (0,1) &~~~~\mbox{ since }~\coth t >1 ~\forall\; t >0.\;
\end{aligned}
\end{equation*}
Now let $\alpha(\kappa,a,s):= \arccos \left(f(\kappa,a,s)\right)$. Then $\alpha(\kappa,a,s) \in (0, \alpha_{\kappa,a})$ $\forall \;\kappa \in \{-1,0,1\}$. Therefore, by Proposition 2.6, there exists an isosceles triangle $T_{\kappa,a,s}$ in $M_\kappa^2$ with base $a$ and base angles $\alpha_{\kappa,a,s}$. Further, equal sides of $T_{\kappa,a,s}$ are $s-a/2$. \hfill $\blacksquare$\\[3mm]
\noindent \textbf{Proposition 2.8}\quad{\it Among all triangles in} $M_\kappa^2$ {\it with base} $a$ {\it and perimeter} $2s_0$ ($s_0 < \pi$ {\it if} $\kappa =1$), 
{\it the isosceles triangle has 
maximum area.} \\[1mm]
\noindent {\it Proof.}\quad By triangle inequality it follows that $s_0>a$. By Proposition 2.7, there exists an isosceles triangle $T_0$ with base $a$ and equal sides $s_0- \dfrac{a}{2}$. Let  $T$ be any triangle in $ M_\kappa^{2} $ with sides $ a, b, c  $ such that $ a+b+c = 2s_0 $. Let $ A,A_0 $ denote the areas of triangles  $ T, T_0 $ respectively. Then, by (2) we get, 
\begin{equation*}
\left.\begin{aligned}
 & T_{|\kappa|}\left( \dfrac{A}{4}\right)  = \sqrt{T_\kappa\left(\dfrac{s_0}{2}\right)\,T_\kappa\left({\dfrac{s_0-a}{2}}\right)\, T_{\kappa}\left({\dfrac{s_0-b}{2}}\right)\,T_\kappa\left({\dfrac{s_0-c}{2}}\right)}~~~~\mbox{and,}\\[1mm]
&T_{|\kappa|}\left( \dfrac{A_0}{4} \right) = \sqrt{T_\kappa\left(\dfrac{s_0}{2}\right)\,T_\kappa\left({\dfrac{s_0-a}{2}}\right)\, T_\kappa\left(\dfrac{a}{4}\right)\,T_\kappa\left(\dfrac{a}{4}\right)}. \end{aligned} \right\}\eqno(11) \end{equation*}
We show that $ A \leq A_0 $ :\quad Note that $\dfrac{A}{4}, \dfrac{A_0}{4} \in I_\kappa$, where $$I_\kappa :=\begin{cases}
(0,\frac{\pi}{4}) & \text{if $\kappa =-1$}\\
(0,\infty) & \text{if $\kappa =0$}\\
(0, \frac{\pi}{2}) & \text{if $\kappa =1$}
\end{cases}$$
and $ T_{|\kappa|} $ is increasing on $ I_\kappa$. Hence $ A \leq A_0 $ if and only if $ T_{|\kappa|} \left(\dfrac{A}{4}\right) \leq T_{|\kappa|}\left(\dfrac{A_0}{4}\right)$. Then by (11) it is enough to verify that $$T_\kappa\left({\dfrac{s_0-b}{2}}\right)\, T_\kappa\left({\dfrac{s_0-c}{2}}\right) \leq T_\kappa^2\left(\dfrac{a}{4}\right).\eqno(12)$$
\noindent Case (i) \underline{$\kappa =0$} :
\begin{eqnarray*}
\mbox{LHS of (12)}&=& 
 \left({\dfrac{s_0-b}{2}}\right)\;\left({\dfrac{s_0-c}{2}}\right)\\
&=&\left({\dfrac{a+c-b}{4}}\right)\;\left({\dfrac{a+b-c}{4}}\right)
= \dfrac{a^2-(b-c)^2}{16}\\
&\leq& \dfrac{a^2}{16}= \left({\dfrac{a}{4}}\right)^2=T_0^2\left({\dfrac{a}{4}}\right)\\[1mm]
&=&\mbox{RHS of (12)}.
\end{eqnarray*}
\noindent Case (ii) \underline{$\kappa \neq 0$} :
\begin{eqnarray*} 
\mbox{LHS of (12)}
&=& 
\dfrac{S_\kappa \left(\dfrac{s_0-b}{2}\right)\, S_\kappa \left(\dfrac{s_0-c}{2}\right)}{C_\kappa \left(\dfrac{s_0-b}{2}\right)\,C_\kappa \left(\dfrac{s_0-c}{2}\right)} 
= -\kappa\; \dfrac{C_\kappa\left(\dfrac{2s_0-b-c}{2}\right)-C_\kappa \left(\dfrac{c-b}{2}\right)}{C_\kappa \left( \dfrac{2s_0-b-c}{2} \right) + C_\kappa \left(\dfrac{c-b}{2}\right)}\\
& &~~~~~~~~~~~~~~~~~~~~~~~~~~~~~~~~~~~~~~~~~~~~~ ~~~\mbox{(\,by (A-11) and (A-12)\,)}\\[1mm]
&=& -\kappa\; \dfrac{C_\kappa\left(\dfrac{a}{2}\right)-C_\kappa \left(\dfrac{c-b}{2}\right)}{C_\kappa \left(\dfrac{a}{2}\right)+ C_\kappa \left(\dfrac{c-b}{2}\right)}.
\end{eqnarray*}
If $b=c=s_0-\dfrac{a}{2}$ then $$\mbox{RHS of (12)}=T_\kappa^2\left(\dfrac{a}{4}\right)
=-\kappa\; \dfrac{C_\kappa\left({\dfrac{a}{2}}\right)-1}{C_\kappa\left({\dfrac{a}{2}}\right)+ 1}.$$ Since $C_\kappa(\theta)\in \begin{cases} [-1,1] & \text{if $\kappa=1$}\\
[1, \infty) & \text{if $\kappa= -1$}\end{cases} $ we get $$\mbox{LHS of (12)}=-\kappa\; \dfrac{C_\kappa\left({\dfrac{a}{2}}\right)-C_\kappa\left({\dfrac{c-b}{2}}\right)}{C_\kappa\left({\dfrac{a}{2}}\right)+ C_\kappa \left({\dfrac{c-b}{2}}\right)} \leq -\kappa\; \dfrac{C_\kappa\left({\dfrac{a}{2}}\right)-1}{C_\kappa\left({\dfrac{a}{2}}\right)+ 1}=\mbox{RHS of (12)}.\eqno{\blacksquare}$$
\noindent \textbf{Theorem 2.9}\quad{\it The following are equivalent for a polygon} $\wp$ {\it in} $M_\kappa^2$ :\\
~\textbf{(i)} \quad $\wp$ {\it is convex.}\\
\textbf{(ii)}\quad $\wp$  {\it is intersection of finitely many closed half-spaces.}\\
\textbf{(iii)}~\,{\it The angle at each vertex of} $\wp$ {\it lies in} $(0, \pi)$.\\[1mm]
\noindent {\it Proof.}\quad $\textbf{(i)} \implies \textbf{(ii)}$ :
Fix $x_0 \in \mbox{interior of}~ \wp$. Let $n$ be the number of vertices of $\wp$. Let $S_1^+, \ldots, S_n^+$ be the closed half-spaces containing $x_0$ corresponding to the boundary geodesic segments $\gamma_1, \ldots, \gamma_n$ of $\partial \wp$ respectively. 
 Then we show that 
$\wp =  S_1^+ \cap \cdots \cap S_n^+$ :\\
$\boldsymbol{\wp \subseteq  
S_1^+ \cap \cdots \cap S_n^+}$ : \quad If not, $\exists \; j \in \{ 1, \ldots , n\}$
  such that $\wp \not \subset S_j^+$. So, $\exists \, y_0 \in$ interior $\wp$ such that $y_0 \not \in S_j^+ $.
  We can assume that $x_0 
\not \in 
\partial S_j^+$. Then by convexity of $\wp$, the convex hull of
  $\{ \gamma_j, \, x_0,\, y_0 \}\subset \wp$, and hence an open set of $M_\kappa^2$ containing midpoint of $\gamma_j$ is also contained in $\wp$.
  This contradicts that $\gamma_j \subseteq 
\partial \wp$. Thus $\wp \subseteq 
S_1^+ \cap \cdots \cap S_n^+$.\\
$ \boldsymbol{S_1^+ \cap \cdots \cap  S_n^+ \subseteq 
\wp} $ : \quad If not, $\exists \,y_0 \in S_1^+ \cap \cdots \cap S_n^+$ such that $y_0 \not \in \wp$\,. Consider $\gamma := [y_0, x_0]$. As $x_0 \in \wp$ and $y_0 \not \in \wp$, $\exists$ a point $z \in \partial \wp \cap \gamma$ such that $[z, y_0]$ intersects $\wp$ only at $z$. Let $i \, \in \{ 1, \ldots, n\}$ be such that $z \in \gamma_i$. Then $y_0 \in M_\kappa^2 \setminus S_i^+$, which gives a contradiction.\\
\noindent $\textbf{(ii)} \implies \textbf{(iii)}$ :
The polygon $\wp$ being an intersection of finitely many closed half-spaces, is convex. Hence at any vertex of $\wp$,
 the angle of the polygon is less than $\pi$.\\
\noindent $\textbf{(iii)} \implies \textbf{(i)}$ :
By $\textbf{(iii)}$, the polygon $\wp$ is locally convex. As $\wp$ is connected, $\wp$ is then a convex polygon.\hfill $\blacksquare$\\[1.5mm]
\noindent The following result follows by Theorem 2.1.\\[0.5mm]
\noindent \textbf{Proposition 2.10}\quad{\it Let} $\wp$ {\it be a convex polygon with} $n$ {\it sides. Let}$~\,\theta_1, \ldots, \theta_n$
{\it be the angles of} $\wp$ {\it at its vertices. Then } $area(\wp) = \kappa \; \left\{ \left(\sum_{i=1}^n \theta_i \right) - (n-2) \pi\right\}$.\\[1.5mm]
\noindent \textbf{Lemma 2.11}\quad{\it The perimeter of any convex} $n${\it -gon in} $M_1^2$ {\it is strictly less than} $2 \pi$.\\[1mm]
\noindent {\it Proof.}\quad Let $\wp$ be a convex $n$-gon with vertices $P_1, \ldots, P_n$ arranged in a cyclic order.
Put $P_{n+1}:= P_1$. Let $a_i$ be the arc-length of the geodesic segment $[P_i, P_{i+1}] \, \forall \, 1 \leq i \leq n$.
 As $\wp$ is a proper polygon, $\{ \underline{0}=:(0,0,0), \, P_i, \, P_{i+1}\}$ determine a plane $H_i$ in $\E^3$ for each
 $i \in \{ 1, \ldots, n \}$. Then $\E^3 \setminus H_i$ has two connected components. We call these components having $H_i$ as common boundary as open half-spaces. Let $H_i^+$  denote the closed half-space in $\E^3$ having $H_i$ as its boundary such that
  $\wp \subset H_i^+$. Then $X:= \cap_{i=1}^n H_i^+$ is a solid cone in $\E^3$ with $\underline{0}$ as its vertex.\\
\par The plane $H$ containing points $ P_1, \, P_2, \, P_3$ intersects $X$
transversely, and $\wp_1 := X \cap H$ is a convex plane-polygon with $n$ sides.  Let $Q_1, \ldots, Q_n$ be the vertices of $\wp_1$ which occur in a cyclic order. Consider the `truncated solid cone' $X_1$ with vertices
$\underline{0}\,, Q_1, \ldots, Q_n$, whose boundary consists of polygon $\wp_1$ and plane-triangles $ \{ \bigtriangleup 
(\,\underline{0}\,,\, Q_i, \,  Q_{i+1}) \}_{1 \leq i \leq n}$.
 (Here, $Q_{n+1}:= Q_1$ and for $1 \leq i \leq n$, $\bigtriangleup 
(\,\underline{0}\,,\, Q_i, \,  Q_{i+1})$ denotes plane-triangle determined by vertices $\underline{0}\,,\, Q_i \;\& \; Q_{i+1}$).
\vspace{1.3cm}
\begin{center}
\psset{unit=1mm}
\begin{pspicture}(-7,-5)(20,15)
\rput(13,1)
{
\includegraphics[height=4cm]{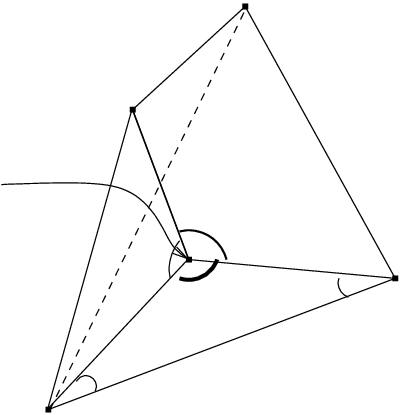}}
\rput(37,-6){$Q_{i+1}$}
\rput(-9.5,3.5){$Q_i$}
\rput(23,-6.8){$\mu_i$}
\rput(16,0.5){$\beta_i$}
\rput(6.5,-4){$\mu_{i-1}$}
\rput(13.2,-8){$\gamma_i$}
\rput(12.5,-17){\large $X_1$}
\rput(11.5,-24){\large Fig. 2}
\rput(1,12){$Q_{i-1}$}
\rput(4.5,-14){ $a_i$}
\rput(-4.5,-20){ $\underline{0}$}
\end{pspicture}
\end{center}
\vspace{2.5cm}
\par Clearly, the face-angles of the polyhedra $X_1$ at the vertex $\underline{0}$ are $a_1, \ldots, a_n$. Thus
$\sum_{i=1}^n a_i$ is the sum of the face angles of $X_1$ at $\underline{0}$. For each $1 \leq i \leq n$, let $\gamma_i$ and $\mu_i$
 be the angles of the plane-triangle $\bigtriangleup 
(\, \underline{0}\,, \, Q_i,\, Q_{i+1} )$ at the vertices $Q_i, \, Q_{i+1}$ respectively. Let $\beta_i$
 be the angle of the polygon $\wp_1$ at vertex $Q_i \, \forall \, 1 \leq i \leq n$. Note that
 $\gamma_i,\, \mu_i,\, \beta_i \in (0, \pi ) \, \forall \, i = 1, \ldots , n$.
 Consider a sphere $S$ with center $Q_i$ having sufficiently small radius $r>0$. Then $S \cap X_1$ is a triangle in $S$  whose sides are of length $r \, \gamma_i,\, r \, \mu_i ~ \& ~ r \, \beta_i $.
 Thus strict triangle inequality holds and we get $\beta_i < \gamma_i + \mu_{i-1}, \, \forall \, i \in \{ 1, \ldots, n\}\,(\mu_0 := \mu_n)$.
$$\therefore~ n \pi  =  \sum_{i=1}^n (a_i + \mu_i + \gamma_i) 
 = \sum_{i=1}^n a_i + \sum_{i=1}^n (\mu_{i-1} + \gamma_i)  > \sum_{i=1}^n a_i + \sum_{i=1}^n \beta_i  = \sum_{i=1}^n a_i + (n-2) \pi.$$
Hence, perimeter of $\wp = \sum_{i=1}^n a_i < 2 \pi$. \hfill $\blacksquare$\\

\noindent \textbf{Lemma of Cauchy (2.12)}\quad {\it Let} $\wp$ {\it and} $\bar{\wp}$ {\it be two convex} $n$-{\it gons in} $M_\kappa^2$ {\it with respective
 vertices} $\{P_i \}_{i= 1, \ldots, n}$ {\it and} $\{ \bar{P_i} \}_{i= 1, \ldots, n}$ {\it occurring in a cyclic order.
  Let} $a_i := d_\kappa(P_i,P_{i+1}) \;\& \;  \bar{a_i}  := d_\kappa(\bar{P_i},\bar{P}_{i+1})\; (1 \leq i \leq n-1)$ {\it denote the lengths of} $(n-1)$ {\it sides of} $\wp$ {\it and} $\bar{\wp}$ {\it respectively. For} $2 \leq i \leq n-1$, {\it let} $\alpha_i$\,({\it resp.} $\bar{\alpha_i}$) {\it denote the angle of} $\wp$\,({\it resp.} $\bar{\wp}$) {\it at vertex} $P_i$\,({\it resp.} $\bar{P_i}$) {\it of} $\wp$\,({\it resp.} $\bar{\wp}$). {\it Let} $a_n$\,({\it resp.} $\bar{a_n}$) {\it be the length of `remaining' side of} $\wp$\,({\it resp.} $\bar{\wp}$). {\it If} $a_i = \bar{a_i}$ {\it for all} $i= 1, \ldots, n-1$ {\it and} $\alpha_i \leq \bar{\alpha_i}$ {\it for all} $i= 2, \ldots, n-1$
 {\it then} $a_n \leq \bar{a_n}$ {\it holds. If in addition, there exists} $i \in \{ 2, \ldots, n-1\}$ {\it with} $\alpha_i < \bar{\alpha_i}$
 {\it then} $a_n < \bar{a_n}$.
\\[1cm]
\begin{center}
\psset{unit=1mm}
\begin{pspicture}(-7,-5)(20,15)
\rput(13,1)
{
\includegraphics[height=4cm]{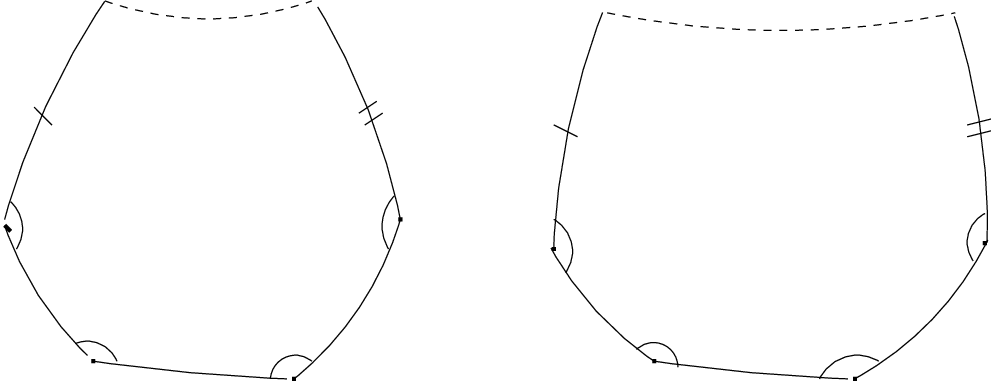}}
\rput(-32,21){$P_1$}
\rput(-18,21){$a_n$}
\rput(-3,21){$P_n$}
\rput(21,21){$\bar{P_1}$}
\rput(41,20){$\bar{a_n}$}
\rput(63,21){$\bar{P_n}$}
\rput(-36,12){$a_1$}
\rput(4,12){$a_{n-1}$}
\rput(18.5,11){$a_1$}
\rput(67.3,11){$a_{n-1}$}
\rput(-19.3,-23){\large $\wp$}
\rput(-3.5,-2.5){$\alpha_{n-1}$}
\rput(7.7,-2.5){$P_{n-1}$}
\rput(15.7,-6){$\bar{P_2}$}
\rput(24,-5){$\bar{\alpha_2}$}
\rput(56.5,-3){$\bar{\alpha}_{n-1}$}
\rput(69,-4){$\bar{P}_{n-1}$}
\rput(39,-23){\large $\tilde{\wp}$}
\rput(-42,-3){ $P_2$}
\rput(-34,-3){ $\alpha_2$}
\rput(15,-24){\large Fig. 3}
\end{pspicture}
\end{center}

\vspace{2.5cm}
\noindent {\it Proof.}\quad By induction on $n$.\\
\noindent \textbf{Claim 1 :}\quad The Lemma of Cauchy is true for $n=3$.\\
From the Law of Cosine for a triangle,
 \begin{center}
\psset{unit=1mm}
\begin{pspicture}(-7,-5)(20,15)
\rput(13,1)
{
\includegraphics[height=2cm]{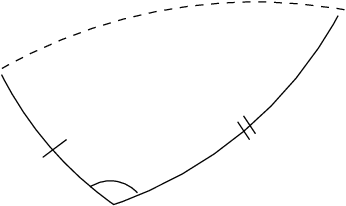}}
\rput(13,12){$a_3$}
\rput(-2,-2.7){$a_1$}
\rput(9,-4.5){$\alpha_2$}
\rput(24,-1.5){$a_2$}
\rput(8,-13){\large Fig. 4}
\end{pspicture}
\end{center}
\vspace{1cm}
$$C_\kappa (a_3) = C_\kappa (a_1) \, C_\kappa (a_2) + \kappa\,S_\kappa (a_1) \, S_\kappa (a_2) \, \cos \alpha_2 ~~~~(\kappa \neq 0),$$
$$ a_3^2 =  a_1^2 +  a_2^2 -2\, a_1 \,  a_2 \, \cos \alpha_2 ~~~~(\kappa = 0).$$Hence it is clear that the side $a_3$ of a triangle
in $M_\kappa^2$ with fixed sides $a_1$ and $a_2$ is a strictly increasing function of $\alpha_2$. This proves Claim 1.\\

Now we assume that the Lemma of Cauchy holds for $n-1\, (n \geq 4)$. Let $\wp$ and $\bar{\wp}$ be two $n$-gons as in the Lemma of Cauchy.\\[1mm]
\noindent \textbf{Claim 2 :}\quad If $\alpha_i = \bar{\alpha_i}$ for some $i \in \{ 2, \ldots, n-1 \}$, then $a_n \leq \bar{a_n}$. Further, if  $\alpha_j < \bar{\alpha_j}$ for some $j \in \{ 2, \ldots, n-1 \} \setminus \{ i\}$, then $a_n < \bar{a_n}$.\\[1mm]
\vspace{1cm}
\begin{center}
\psset{unit=1mm}
\begin{pspicture}(-7,-5)(20,15)
\rput(13,1)
{\includegraphics[height=3.5cm]{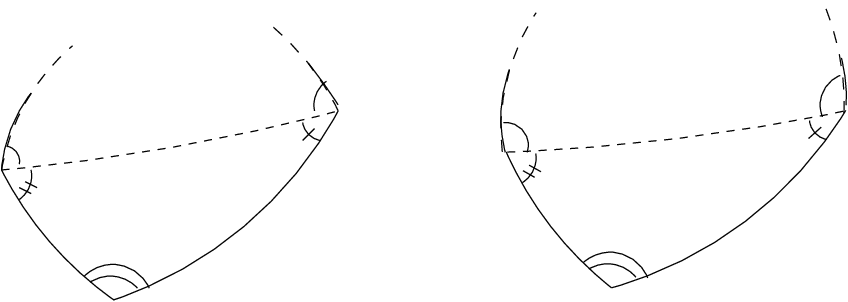}
}
 \rput(-40,-2){$P_{i-1}$}
 \rput(-35,-9.5){$a_{i-1}$}
 \rput(-4.5,8.5){$\alpha'_{i+1}$}
 \rput(-29,3){$\alpha'_{i-1}$}
 \rput(-21,-10){$\alpha_i$}
 \rput(-18,10){\large $\wp'$}
\rput(29.5,5.5){$\bar{\alpha'}_{i-1}$}
\rput(39,12){\large $\bar{\wp'}$}
\rput(54,10){$\bar{\alpha'}_{i+1}$}
\rput(66.5,5.5){$\bar{P}_{i+1}$}
\rput(-2,-6.7){$a_i$}
\rput(23.5,-8.5){$a_{i-1}$}
\rput(36,-9){$\alpha_i$}
\rput(54,-7.7){$a_i$}
\rput(8,6){$P_{i+1}$}
\rput(19,-1){$\bar{P}_{i-1}$}
\rput(-16,-4.7){\large $T$}
\rput(42,-3){\large $\bar{T}$}
\rput(-23,-19.7){$P_i$}
\rput(34.5,-18.9){$\bar{P_i}$}
\rput(-18,-25){\large $\wp$}
\rput(40,-25){\large $\bar{\wp}$}
\rput(10,-28){\large Fig. 5}
\end{pspicture}
\end{center}
\vspace{2.5cm}

\quad Let $\gamma_i$ and $\bar{\gamma_i}$ denote the geodesic segments 
 $[P_{i-1}, \, P_{i+1}]$ and $[\bar{P}_{i-1}, \, \bar{P}_{i+1}]$ respectively. Since $\wp$ and $\bar{\wp}$  are convex we
 obtain two convex $(n-1)$-gons $\wp'$ and $\bar{\wp}'$ and two 
 triangles $T$ and $\bar{T}$ as shown in the figure above. Note
that $T$ and $\bar{T}$ are congruent triangles. In particular, 
the angles of $T, \, \bar{T}$ at $P_{i-1}$ and $\bar{P}_{i-1}$(resp. at $P_{i+1}$ and $\bar{P}_{i+1}$) are equal. This implies that 
 $\alpha'_{i-1} \leq \bar{\alpha}'_{i-1}$ and  $\alpha'_{i+1} \leq \bar{\alpha}'_{i+1}$. Here, $\alpha'_{i-1}$, $\alpha'_{i+1}$ are angles of $\wp'$
at vertices $P_{i-1}, \, P_{i+1}$ respectively. Similarly 
$\bar{\alpha}'_{i-1}$, $\bar{\alpha}'_{i+1}$ are defined. 
 Thus $\wp'$ and $\bar{\wp}'$
satisfy the assumption of the Lemma of Cauchy, and Claim 2 follows by induction assumption.\\

By Claim 2, we can now assume that $\wp$ and $\bar{\wp}$ are two convex $n$-gons as in the Lemma of Cauchy which further satisfy $\alpha_i < \bar{\alpha_i}\; \forall \, i \in \{ 2, \ldots, n-1 \}$. We show that $a_n < \bar{a_n}$ :
Increase the angle $\alpha_{n-1}$
of $\wp$ at $P_{n-1}$ 
 until it becomes equal to $\bar{\alpha}_{n-1}$, while maintaining the $(n-1)
$ sides constant. This way, we obtain a new polygon $\wp'$ with vertices $P_1, \, \ldots, \, P_{n-1}, \, P'_n$,\,
side-lengths $a_1, \, \ldots, \, a_{n-1}, \, a'_n := d_\kappa(P'_n,P_1)$ 
and angles at vertices $P_2, \, \ldots ,\, P_{n-2}, \, P_{n-1}$
 equal to $\alpha_2, \, \ldots ,\, \alpha_{n-2}, \, \bar{\alpha}_{n-1}$ respectively.\\

\noindent \textbf{Case (i)}\quad $\wp'$ is convex :\\
Join $P_1, \, P_{n-1}$ by geodesic segment $\gamma$ (say). Since $\wp'$ is convex, $\gamma \subset \wp'$ and $\gamma$
 divides $\wp'$ into two convex proper polygons. Apply Claim 1 to the two triangles $[P_1,P_{n-1},P_n]$ and
 $[P_1,P_{n-1},P'_n]$, whence $a_n = d_\kappa(P_1,P_n) \leq d_\kappa(P_1,P'_n) =a'_n$, and 
$a_n < a'_n$ since 
$\alpha_{n-1}< \bar{\alpha}_{n-1}$.\\
\vspace{1.1cm}
\begin{center}
\psset{unit=1mm}
\begin{pspicture}(-7,-5)(20,15)
\rput(13,1)
{\includegraphics[height=4cm]{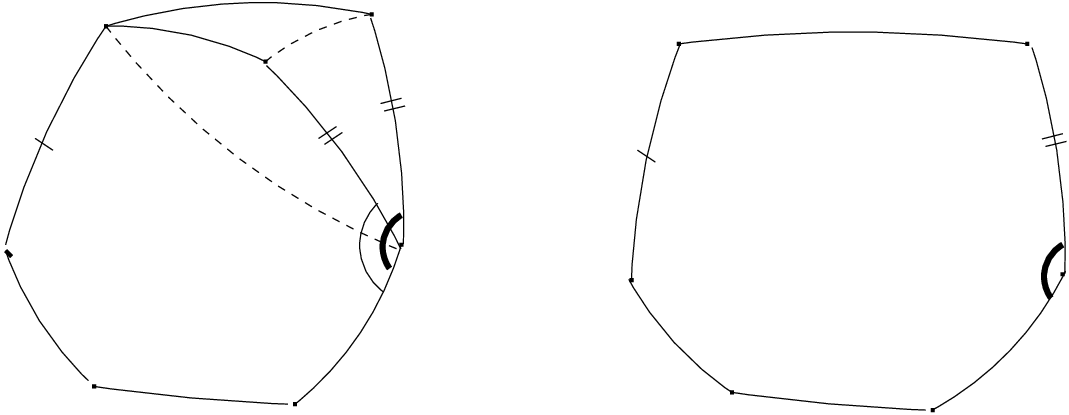}
}
\rput(-31,21){$P_1$}
\rput(-18,23.5){$a'_n$}
\rput(0,21){$P_n'$}
\rput(-20,16){$a_n$}
\rput(-14,12.5){$P_n$}
\rput(23.5,18.5){$\bar{P_1}$}
\rput(42,20){$\bar{a_n}$}
\rput(63.5,18.5){$\bar{P_n}$}
\rput(-19.5,-23.5){ $\large{\wp'}$ (convex)}
\rput(-8.5,-2.5){$\alpha_{n-1}$}
\rput(5.6,-2.5){$P_{n-1}$}
\rput(18.5,-6){$\bar{P_2}$}
\rput(-19,5){$\gamma$}
\rput(57,-3.5){$\bar{\alpha}_{n-1}$}
\rput(69,-5.5){$\bar{P}_{n-1}$}
\rput(41,-23){\large $\bar{\wp}$}
\rput(-42,-3){ $P_2$}
\rput(13,-23.5){\large Fig. 6}
\end{pspicture}
\end{center}

\vspace{2.5cm}
\par We then apply the induction assumption and Claim 2 to the $n$-gons $\wp'$ and $\bar{\wp}$ which
have the same angles at $P_{n-1}$ and $\bar{P}_{n-1}$, 
 and obtain
$$a'_n = d_\kappa(P_1,P'_n) \leq d_\kappa(\bar{P_1},\bar{P_n}) =\bar{a_n}.$$ Thus $a_n <a'_n \leq \bar{a_n}$ and we have concluded the proof for case (i).\\

\noindent \textbf{Case (ii)}\quad $\wp'$ is not convex :\\
In this case, as we increase $\alpha_{n-1}$ by rotating side $[P_{n-1},P_n]$ around $P_{n-1}$, there exists a smallest value
$\alpha'_{n-1}$ of the angle for which $\wp'$ ceases to be convex. This value lies between
$\alpha_{n-1}$ and $\bar{\alpha}_{n-1}$.\\[3mm]
\begin{center}
\psset{unit=1mm}
\begin{pspicture}(2,-5)(20,15)
\rput(13,1)
{\includegraphics[height=3cm]{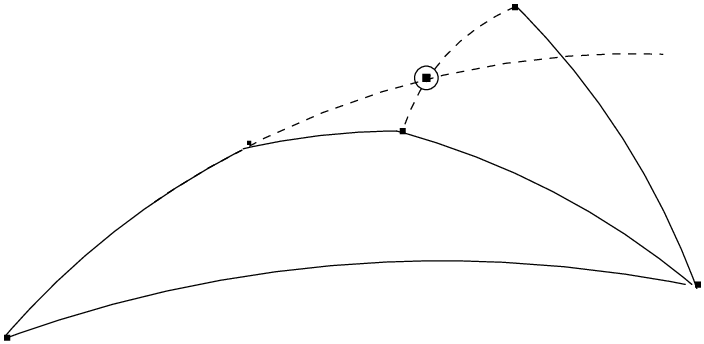}
}
\rput(0.5,5){$P_1$}
\rput(28,18.5){$P_4'$}
\rput(17,12.7){$P_4''$}
\rput(16.5,1){$P_4$}
\rput(46,-10){$P_3$}
\rput(-21,-15){ $P_2$}
\rput(18,-14){$\large \wp_4'$ (not convex)}
\rput(19,-19.5){\large Fig. 7}
\end{pspicture}
\end{center}

\vspace{2cm}
Let $P''_n$  be the point thus obtained. By construction, $P''_n$ belongs to the
{\it line} 
determined by $P_2$ and $P_1$. We have $$a''_n := d_\kappa(P_1,P''_n) = d_\kappa(P_2,P''_n) - d_\kappa(P_1,P_2) = d_\kappa(P_2,P''_n) - a_1  \eqno{(13)}$$\\[1cm]
\begin{center}
\psset{unit=1mm}
\begin{pspicture}(-7,-5)(20,15)
\rput(13,1)
{\includegraphics[height=5cm]{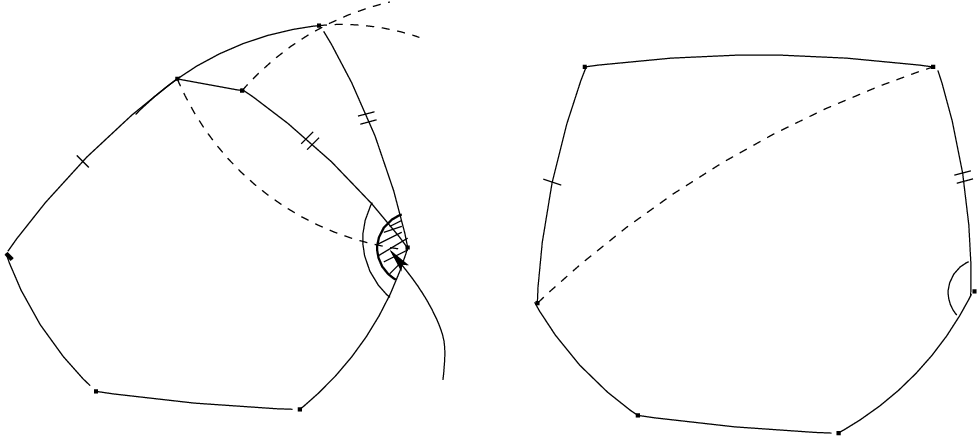}
}
\rput(-26,19){$P_1$}
\rput(-15,24.5){$a''_n$}
\rput(-7,26.1){$P''_n$}
\rput(-17.2,17.7){$a_n$}
\rput(-15.5,13){$P_n$}
\rput(40,4){$\bar{\gamma}$}
\rput(21,20){$\bar{P_1}$}
\rput(41,22){$\bar{a_n}$}
\rput(65.5,21){$\bar{P_n}$}
\rput(-21.3,-12){\large $\wp$}
\rput(7,-19){$\alpha'_{n-1}$}
\rput(-6,-4.5){$\alpha_{n-1}$}
\rput(8.3,-2.5){$P_{n-1}$}
\rput(15.7,-9){$\bar{P_2}$}
\rput(73,-8){$\bar{P}_{n-1}$}
\rput(42,-14){\large $\bar{\wp}$}
\rput(-7,-23.5){$P_{n-2}$}
\rput(-46,-3){ $P_2$}
\rput(9,-27){\large Fig. 8}
\end{pspicture}
\end{center}

\vspace{3cm}

\noindent Applying triangle inequality to the 
triangle $[\bar{P_1},\bar{P_2},\bar{P_n}]$ we get
\begin{eqnarray*}
\qquad \qquad \qquad \quad ~\bar{a_n} = d_\kappa(\bar{P_1},\bar{P_n}) & \geq &  d_\kappa(\bar{P_2},\bar{P_n}) - d_\kappa(\bar{P_1},\bar{P_2}) \\
& = &  d_\kappa(\bar{P_2},\bar{P_n}) - d_\kappa(P_1,P_2) \\
& = &  d_\kappa(\bar{P_2},\bar{P_n}) - a_1. \hfill{\qquad ~~\qquad \qquad \qquad \quad}(14)
\end{eqnarray*}
Now we can apply induction assumption to the convex $(n-1)$-gons $[\bar{P_2},\ldots,\bar{P_n}]$ and $[P_2,P_3, \ldots, P_{n-1},P''_n]$ to get
$$d_\kappa(\bar{P_2},\bar{P_n}) \geq d_\kappa(P_2, P''_n). \eqno{(15)} $$ 
Finally, applying Claim 1 to the triangles $[P_1,P_n,P_{n-1}]$ and $[P_1,P''_n,P_{n-1}]$ we get
$$ a''_n 
>                             a_n \eqno{(16)}$$
Thus,
\begin{alignat*}{2}
\bar{a_n} 
& \geq   d_\kappa(\bar{P_2},\bar{P_n}) - a_1    & \qquad & [\,\mbox{by} ~(14)\,]\\
& \geq   d_\kappa(P_2,P_n'') - a_1 & \qquad & [\, \mbox{by} ~(15)\, ]\\  
& =   a''_n 
& \qquad & [\, \mbox{by} ~ (13)\, ]\\
& >  
a_n 
& \qquad & [\, \mbox{by} ~(16)\,].
\end{alignat*}
This proves the Lemma of Cauchy. \hfill $\blacksquare$\\[3mm]
\noindent \textbf{Lemma 2.13}\quad{\it Among all convex} $n${\it -gons in} $M_\kappa^2$ {\it whose all sides but
one are given in length -- say}
 $a_1, \, \ldots, \, a_{n-1}$
 {\it --} ({\it with} $a_1 +  \cdots + a_{n-1} < \pi$ {\it if} $\kappa=1$), {\it area maximizer
 is the convex} $n${\it -gon whose
 vertices lie on a circle having its center at the midpoint of the remaining side.}\\[1mm]
\noindent {\it Proof.}\quad When $n=3$, this result is proved in Proposition 2.5. Here we consider $n \geq 4$.
Let $\U$ denote the family of all convex $n$-gons in $M_\kappa^2$ whose all sides but one are
 $a_1, \, \ldots, \, a_{n-1}$ (with $a_1 +  \cdots + a_{n-1} < \pi$ if $\kappa=1$). Let $r_0:= a_1 +  \cdots + a_{n-1}$.
\noindent {\it Existence} \quad Upto congruence all polygons in $\U$ lie inside $B_\kappa(p,r_0)$ where $p \in M_\kappa^2$. Since $\overline{B_\kappa(p,r_0)}$ is compact in $(M_\kappa^2, d_\kappa)$ and number of vertices is $n$ for all polygons in $\U$,
 there exists an `area maximizer' $\wp$ in $\U$.\\

 Let
 $A_1, \, \ldots, \, A_n \,(A_{n+1}:= A_1)$ be the vertices of $\wp$ which occur in a cyclic order
 and such that $d_\kappa(A_i,A_{i+1})= a_i \,\forall \, i= 1,\ldots, \, n-1$.
Put $r:= d_\kappa(A_1,A_n)/2$ and $O:=$ mid-point of geodesic segment $[A_1,A_n]$. We show that
$d_\kappa(O,A_i)= r \; \forall \; i = 1, \, \ldots, \, n$ :
\par Suppose
$d_\kappa(O,A_i)\neq r $ for some $i \in \{2, \, \ldots, \, n-1 \}$. Put $a:= d_\kappa(A_1,A_i)$ and
$b:= d_\kappa(A_n,A_i)$. Since $d_\kappa$ is a metric on $M_\kappa^2$ we get $a+b \leq
 a_1+a_2+  \cdots  a_{n-1}$ ($< \pi$ if $\kappa =1$), and by assumption $A_i$ does not lie on the circle of radius
 $ r $ and center $O$. By Proposition 2.5, there exists a triangle $[A'_1,A_i,A'_n]$ such that
$d_\kappa(A'_1,A_i)=a, ~ d_\kappa(A'_n,A_i)=b$ and $$area(
[A'_1,A_i,A'_n]
) > area([A_1,A_i,A_n]) \eqno{(17)}$$
Further we can assume that the angles of
$[A'_1,A_i,A'_n]$ at vertices
$A'_1,\, A_i,\, A'_n$ are close to the angles of
$[A_1,A_i,A_n]$ at vertices
$A_1,\, A_i,\, A_n$ respectively.
\par Let $T\, '$ be the triangle $[A'_1,A_i,A'_n]$. Let $S_1$ be the closed half-space of $M_\kappa^2$ 
containing $A'_n$
and having the {\it line} 
containing $[A'_1,A_i]$ as its boundary. Let $S_2$ denote the other closed half-space. Consider the polygon $\wp'_1 \subset S_2$ with vertices $A'_1, A'_2, \ldots, A'_i$ with $A'_i = A_i$, occurring in a cyclic order such that
$\wp'_1$ is congruent to the convex polygon $[A_1, A_2, \ldots, A_i]$.  Let $S_3$ be the closed half-space of $M_\kappa^2$ containing $A'_1$
and having the {\it line} containing $[A_i,A'_n]$ as its boundary. Let $S_4$ denote the other closed half-space. Similarly, consider a polygon
$\wp'_2 \subset S_4$ with vertices $A'_i (= A_i),A'_{i+1}, \ldots, A'_n$ occurring in a cyclic order such that
$\wp'_2$ is congruent to the convex polygon $[A'_i, A_{i+1}, \ldots, A_n]$. Polygons $\wp'_1$ and $\wp'_2$ do not intersect the
interior of $T \, '$. Thus we have constructed a polygon $\wp'$, with $n$ vertices $A'_1,\, A'_2,\, \ldots , \, A'_n$
occurring in a cyclic order and such that $\wp' = \wp'_1 \cup T\, ' \cup \wp'_2$. By $(17)$ and
construction of $\wp'$, $area(\wp') > area(\wp)$. Also, since the angles of $T \, '$
 at vertices $A'_1, \, A_i, \, A'_n$ are {\it sufficiently close} to the angles of $[A_1, A_i,A_n]$ at vertices
$A_1,\,  A_i,\, A_n$ respectively, then the angles of $\wp'$ at the vertices $A'_1, \, A'_i = A_i, \, A'_n$  are strictly less than
$\pi$. By Theorem 2.9 it follows that $\wp'$ is a {\it convex polygon}. Thus $\wp' \in \U$ and $area(\wp')> area(\wp)$,
which contradicts the fact that $\wp$ is an `area maximizer' in  $\U$. We conclude that $d_\kappa(O,A_i) =r \; \forall \; i \in \{ 1, \ldots , n\}$.\hfill $\blacksquare$ \\[2mm]
\noindent \textbf{Lemma 2.14}\quad{\it Let} $\C$ {\it be any piecewise smooth closed curve in} $M_1^2$ {\it
whose arc-length is strictly less than} $2 \pi$. {\it Then} $\C$ {\it is contained in an open hemisphere}. \quad (cf. \cite{Mil-Park})\\[2mm]
\noindent \textbf{Definition :}\quad A {\it digon} $D_{x, \alpha}$ ($x \in M_1^2$ and $\alpha \in [0, \pi]$) is a closed region of $M_1^2$ bounded by two half great circles with end points $x, -x$ and forming an angle $\alpha$ at $x$.\\[2mm] \noindent \textbf{Remark :}\quad The area of the digon $$D_{p_0,\alpha}=\left\{(\cos \theta~ \cos \phi,~ \cos \theta ~ \sin \phi,~ \sin \theta) : \theta \in
\left[- \frac{\pi}{2}, \frac{\pi}{2} \right], \phi \in [0, \alpha] \right\} $$
is equal to $2 \alpha$ since $ \displaystyle \int_0^\alpha \int_{-\frac{\pi}{2}}^{\frac{\pi}{2}} \cos \theta \, d \theta \, d \phi = 2 \alpha $. There is an isometry between any two digons with the same angle $\alpha$. Hence for each $\alpha \in [0,\pi]$, area of digon $D_{x,\alpha}$ is $2\,\alpha$ $\forall \;x \in M_1^2$.\\[2mm]
\noindent \textbf{Lemma 2.15}\quad{\it Let} $\C$ {\it be a piecewise smooth closed curve in} $M^2_1$ {\it with arc-length} $2 \pi$. {\it If} $\C$ {\it is not a digon then} $\C $ {\it is contained in an open hemisphere}. \quad (cf. \cite{Mil-Park})
\begin{center}
3. REGULAR POLYGONS IN $M_\kappa^2$ 
\end{center}
 A polygon in $M^2_\kappa$ is said to be {\it equilateral} (resp. {\it equiangular}) if all its sides
have same length\,(resp. if all its angles are equal). A polygon is said to be {\it regular} if it is convex, equilateral and equiangular.
A regular polygon (proper regular polygon if $\kappa=1$) of $n$ sides is called a {\it regular n-gon}.\\[2mm]
\noindent \textbf{Construction of regular polygons in} $\boldsymbol{M_\kappa^2}$ \textbf{:}\quad Fix $r>0$ ($r < \frac{\pi}{2}$ if $\kappa =1$) and $n \geq 3$. Let $p_0 \in M_\kappa^2$ be as in (1). Let $\C_\kappa(p_0,r) $ denote the circle which is the boundary of the disc $B_\kappa(p_0\, , r)$ contained in $M_\kappa^2$. Then $\C_\kappa(p_0,r)$ is nothing but 
a Euclidean circle in the plane $\{\left(x_1,x_2, |\kappa|\, C_\kappa (r)\right)\, |\; x_1, x_2 \in \R\}\subseteq \R^3$ with center $ \underline{c}=C_{\kappa} (r) \, p_0$ and radius $S_{\kappa} (r)$. \\[1mm]
Let $P_1, \ldots ,P_n$ be $n$ points in $\C_\kappa(p_0,r)$
which occur clockwise such that
$\measuredangle \{ P_i - \underline{c} \; , P_{i+1} - \underline{c}\} = \frac{2 \pi}{n}\; \forall \; i = 1,\ldots ,n$ (\,here, $P_{n+1} := P_1$\,). Let $\wp_{n,r}$ denote the convex polygon in $M_\kappa^2$ with $P_1, \ldots ,P_n$ as its vertices.
 By construction, the rotation, $\rho_{\frac{2 \pi}{n}}$ about the oriented axis through $\underline{c}$ normal to the plane of $\C_\kappa(p_0,r)$ is a symmetry of $\wp_{n,r}$. Here the axis is oriented by the vector $(0,0,1)$. 
Thus, $\wp_{n,r}$ is an equilateral, equiangular $n$-gon in $M_\kappa^2$. Any two convex polygons constructed as above are congruent
 to each other for a fixed $n \geq 3$  and fixed $r>0$ ($r < \frac{\pi}{2}$ if $\kappa=1$). \\[2mm]
Let $a$ be the length of a side of $\wp_{n,r}$. Let $Q$ be the midpoint of $[P_1,P_2]$. Then the triangles $[\underline{c},Q, P_1]$ is congruent to the triangle $[\underline{c},Q, P_2]$ and for both these triangles the angle at the vertex $Q$ is $\pi/2$. The Law of Sine (B-3) applied to the triangle $[\underline{c},Q, P_1]$ gives
$$\dfrac{S_\kappa(r)}{\sin \left(\frac{\pi}{2}\right)}= \dfrac{S_\kappa\left(\frac{a}{2}\right)}{\sin \left( \frac{\pi}{n}\right) }.$$
Therefore,  $$a = a(n,r)  = 2\, AS_\kappa \left( S_\kappa (r) \; \sin \left( \frac{\pi}{n}\right) \,\right). \eqno{(18)}$$ 
  Now we compute the angle $\theta = \theta (n,r)$ at vertices of $\wp_{n,r}$. Recall that $n \geq 3$. The Law of Sine (B-3) applied to the triangle $[P_1, P_2, \underline{c}]$ we get 
$$\dfrac{S_\kappa(r)}{\sin \left( \frac{\theta}{2}\right) }= \dfrac{S_\kappa(a)}{\sin \left( \frac{2\,\pi}{n}\right) }.\eqno{(19)}$$
\noindent \underline{$\kappa=0$} \quad From (18) and (19) it follows that 
$$\sin \left( \frac{\theta}{2} \right) =  \dfrac{ 2\, \sin \left( \frac{\pi}{n}\right) \, \cos \left( \frac{\pi}{n}\right) \, S_\kappa(r)}{2\, \sin \left( \frac{\pi}{n}\right) \, S_\kappa(r)}=\cos \left( \dfrac{\pi}{n} \right) =\sin \left(\dfrac{\pi}{2}- \dfrac{\pi}{n}\right) .
$$
Therefore,
  $$ \theta = \theta(n,r) = \left( \dfrac{n-2}{n}\right)\, \pi ~~~~~~~~(\kappa =0).\eqno{(20)} $$\\
\noindent \underline{$\kappa \neq 0$}\quad From (18) and (19) we get,
$$\dfrac{S_\kappa(r)}{\sin \left( \frac{\theta}{2}\right) } = \dfrac{S_\kappa(a)}{\sin \left( \frac{2\,\pi}{n}\right) }
=
\dfrac{2\, S_\kappa\left(\frac{a}{2}\right)\, C_\kappa\left(\frac{a}{2}\right)}{2\,\sin \left( \frac{\pi}{n}\right) \, \cos  \left( \frac{\pi}{n}\right) }
=
 \dfrac{ S_\kappa(r)\, \sin \left( \frac{\pi}{n}\right) \,C_\kappa\left(\frac{a}{2}\right)}{\sin \left( \frac{\pi}{n}\right) \, \cos  \left( \frac{\pi}{n}\right) }.
$$
Thus $$\sin \left( \dfrac{\theta}{2} \right) = \dfrac{ \cos \left(  \frac{\pi}{n}\right) }{C_\kappa\left(\frac{a}{2}\right)}
 ~~~~~~~~(\kappa \neq 0).\eqno{(21)}$$
Therefore, by (21) and (18), 
\begin{eqnarray*}
\cos \left( \dfrac{\theta}{2} \right) 
&=& \dfrac{ \sqrt{C_\kappa^2\left(\frac{a}{2}\right)-\cos^2  \left( \frac{\pi}{n}\right) }}{C_\kappa\left(\frac{a}{2}\right)}
=
\dfrac{ \sqrt{1-\kappa\, S_\kappa^2\left(\frac{a}{2}\right)-\cos^2  \left(\frac{\pi}{n}\right)}}{C_\kappa \left(\dfrac{a}{2}\right)}\\
&=& 
\dfrac{ \sqrt{\sin^2  \left( \frac{\pi}{n}\right) -\kappa\, S_\kappa^2(r)\, \sin^2\left(\frac{\pi}{n}\right)}}{C_\kappa \left(\frac{a}{2}\right)}\\
&=& \dfrac{ \sqrt{\sin^2  \left( \frac{\pi}{n}\right) \left( 1-\kappa\, S_\kappa^2(r)\right)}}{C_\kappa \left(\frac{a}{2}\right)}
= 
\dfrac{\sqrt{ \sin^2  \left( \frac{\pi}{n}\right) \,C_\kappa^2(r)}}{C_\kappa \left(\frac{a}{2}\right)}.
\end{eqnarray*}
Thus $$\cos \left( \dfrac{\theta}{2} \right) = \dfrac{ \sin  \left( \frac{\pi}{n}\right) \,C_\kappa(r)}{C_\kappa \left(\frac{a}{2}\right)} ~~~~~~~~(\kappa \neq 0).\eqno(22)$$
From (21) and (22) we get,
$$\tan \left( \dfrac{\theta}{2}\right)=  \dfrac{ \cos \left( \frac{\pi}{n}\right) }{\,C_\kappa(r)\,\sin \left( \frac{\pi}{n}\right) } ~~~~~~~~(\kappa \neq 0).$$
Therefore $$\theta= \theta(n,r)= 2\, \arctan \left( \dfrac{\cot\left( \dfrac{\pi}{n}\right)}{C_\kappa(r)} \right) ~~~~~~~~(\kappa \neq 0).\eqno{(23)}$$
Let $A$ denote the area of the regular $n$-gon $\wp_{n,r}$. \\
\noindent \underline{$\kappa \neq 0$} : \quad By Proposition 2.10, $A= \kappa \; \left\{
n \,\theta  - (n-2) \pi\right\}$. Therefore by (23),
$$A= A(n,r)=\kappa \; \left\{ 
2\, n \,\arctan \left( \dfrac{\cot\left( \frac{\pi}{n}\right)}{C_\kappa(r)} \right)   - (n-2) \pi\right\} ~~~~~~~~ (\kappa \neq 0).\eqno{(24)}$$
\noindent \underline{$\kappa = 0$} :\quad Let $T$ be the triangle $ [\underline{c}, P_1,P_2]$. Then $ 
area(T)= 
\dfrac{1}{2} \, r^2 \, \sin \left(\dfrac{2\,\pi }{n}\right)$. Hence, $$A= A(n,r)=n\,\, area(T)=  \dfrac{n\,r^2}{2} \sin \left(\dfrac{2\,\pi }{n}\right) ~~~~~~~~(\kappa =0).\eqno{(25)}$$ 
\noindent \textbf{Theorem 3.1}\quad{\it Any regular} $n${\it -gon in} $M_\kappa^2$ {\it is congruent to} $\wp_{n,r}$ {\it for a unique}
$r>0$ ($r< \frac{\pi}{2}$ {\it if} $\kappa =1$).\\[1mm]
\noindent {\it Proof.}\quad Let $\wp'$ be any regular $n$-gon in $M_\kappa^2$. As $\wp'$
is a regular $n$-gon, $n \geq 3$ holds. Let $a'$ be the length of a side of $\wp'$. By Lemma 2.11, $n a' < 2 \pi$ if $\kappa =1$. That is, $$a' \in \begin{cases}
(0, \frac{2\, \pi}{n}) & \text{if $\kappa =1$}\\
(0, \infty) & \text{if $\kappa \neq 1$}. 
  \end{cases}$$ Let
$$ J_\kappa := \begin{cases}
(0,  \pi/2) & \text{if $\kappa =1$}\\
(0, \infty) & \text{if $\kappa \neq 1$} 
  \end{cases} ~~\mbox{ and }~~ J_\kappa^\prime := \begin{cases}(0,  2 \pi/n) & \text{if $\kappa =1$}\\
(0, \infty) & \text{if $\kappa \neq 1$}. 
  \end{cases}$$ 
Consider the function $f: J_\kappa
\rightarrow J_\kappa^\prime
$ defined by
$f(r) := 2\, AS_\kappa \left( \sin\left(\dfrac{\pi}{n} \right)\, S_\kappa(r) \right)
$. 
Then $f: J_\kappa
\rightarrow J_\kappa^\prime
$ is a bijection for a fixed $n \geq 3$. Hence for $a' \in J_\kappa^\prime,\; \exists$ a unique $r \in J_\kappa$ such that $a' = f(r)$. Thus $a' = a(n,r)$
 for a unique $r \in J_\kappa$.
\par Now we prove that $\wp'$ is congruent to $\wp_{n,r}$. Let $P'_1, \ldots , P'_n$ (resp. $P_1, \ldots, P_n$) be the vertices of
$\wp'$\,(resp. of $\wp_{n,r}$) which occur in a cyclic order. Let $\theta \, '$\,(resp. $\theta$) be the angle of $\wp'$\,(resp. of $\wp_{n,r}$) at
 its vertices. If $\theta \, ' < \, \theta$\,(resp. $> \theta$), then by the {\it Lemma of Cauchy}\;(Lemma 2.12),
 $d_\kappa(P'_1,P'_n) < d_\kappa(P_1,P_n)$\;(\,resp. $> d_\kappa(P_1,P_n)$\,) which contradicts that $a' = a(n,r)$. So, $\theta \, ' = \theta$.
 Applying the {\it Lemma of Cauchy} again to the convex polygons $[P_1, P_2. \ldots, P_j] \; \& \;[P'_1, P'_2. \ldots, P'_j]$,
 we get $d_\kappa(P'_1,P'_j) = d_\kappa(P_1,P_j) \, \forall \; j = 2, \ldots , n$. Similarly,
 $d_\kappa(P'_i,P'_j) = d_\kappa(P_i,P_j) \, \forall \; i,j \in \{ 1, \ldots , n\}$. By Proposition 1.1, there exists an isometry $\varphi$
 of $M_\kappa^2$ such that $\varphi(\wp') = \wp_{n,r}$. \hfill $\blacksquare$\\[2mm]
\noindent \textbf{Proposition 3.2}\quad{\it Let} $\wp$ {\it be a regular} $n${\it -gon in} $M_\kappa^2$
{\it having side} $a$, {\it angle} $\theta$ {\it and area} $A$. {\it Then}
$\exists$ a unique $r>0$ ($r <\frac{\pi}{2}$ {\it if} $\kappa =1$) {\it such that}
$\wp$ {\it is inscribed in a circle of radius} $r$. {\it Further, equations}\\
\begin{align*}
\mbox{\textbf{(i)}}  \quad & r = r(n,a) = AS_\kappa \left(  \frac{S_\kappa (a/2)}{\sin (\pi/n)}\right),
\qquad \qquad \qquad \qquad \qquad \qquad \qquad \qquad \quad \,\hfill (26)\\
\mbox{\textbf{(ii)}} \quad & r = r(n,\theta) = AC_\kappa \left(
 \frac{ \cot (\pi/n)}{\tan (\theta/2)}
\right) 
\qquad ~~ (\kappa \neq 0),  \qquad \quad  \qquad \qquad \qquad \qquad  \,\hfill (27)\\
\mbox{\textbf{(iii)}}  \quad & r = r(n,A) = 
\begin{cases}
AC_\kappa \left(\cot \left(\frac{\pi}{n}\right)\,\tan \left( \dfrac{2 \pi - \kappa\, A}{2n}\right)\right) & \text{if $\kappa \neq 0$} \qquad \qquad \quad \,~\hfill(28)\\[3mm]
\sqrt{\dfrac{2\,A}{n\, \sin\left( \frac{2\,\pi}{n}\right) } }\,& \text{if $\kappa = 0$} \qquad \qquad \quad ~\,\hfill(29)
\end{cases}
\end{align*}
{\it hold, and any regular} $n${\it -gon in} $M_\kappa^2$ ($\kappa \neq 0$) {\it is determined (uniquely up to congruence) by any one of three} :
$$a \in
\begin{cases}
(0,\pi)~~~ \text{ {\it if }~ $\kappa =1$}\\
(0, \infty)~~ \text{ {\it if}~ $\kappa =-1$}
\end{cases}, ~~\; \theta \in (0, \pi) 
,~~~A \in
\begin{cases}
(0, 2 \, \pi)~~~~~~~~~~ \text{ {\it if }~ $\kappa =1$}\\
(0, (n-2)\,\pi)~~ \text{ {\it if}~ $\kappa =-1$}
\end{cases}.$$ {\it Further, any regular} $n${\it -gon in} $M_0^2$ {\it is determined (uniquely up to congruence) by any one of two} :
$a,A\in (0,\infty)$.\\[1mm]
\noindent {\it Proof.}\quad By the Theorem 3.1, there exists a unique $r >0$ ($r< \frac{\pi}{2}$ if $\kappa=1$) such that $\wp$ is congruent to $\wp_{n,r}$. Now equation (26), (27), (28) and (29) easily follows from (18), (23), (24) and (25) respectively. 

\par Finally, it can be verified that the functions in \textbf{(i)}, \textbf{(ii)}, \textbf{(iii)} are strictly monotone functions and hence
$\wp$ is determined up to congruence by any one of the entities $\theta$ (if $\kappa \neq 0$), $a$ and $\, A$. \hfill $\blacksquare$\\[3mm]
\noindent \textbf{Remark :} Fix $n \in \N$. When $\kappa =0$ the angle $\theta \in (0, \pi)$ is not enough to determine the regular $n$-gon. 
But for $\kappa \neq 0$ the regular $n$-gons in $M_\kappa^2$ with angle $\theta$ are congruent.\\[3mm]
\noindent \textbf{Corollary 3.3}\quad{\it Let} $(\tilde{\wp_k})_{k \in \N}$ {\it be a sequence of regular polygons in} $M_\kappa^2$ ({\it proper regular polygons if} $\kappa=1$)
{\it such that} $\tilde{\wp_k}$  {\it has} $k$ vertices $\forall \; k$, and $A_k:=(area(\tilde{\wp_k})) \longrightarrow A'$ {\it as} $k \longrightarrow \infty$. {\it For each}
 $k \in \N$, {\it let} $r_k$
{\it be the radius of the circle in which} $\tilde{\wp_k}$ {\it is inscribed. Then}
\begin{align*}
\textbf{(i)} \quad & \lim_{k \longrightarrow \infty} r_k = AS_\kappa \left(\sqrt{A' \, (4 \pi -\kappa\, A')}/(2 \pi)\right)\\
\textbf{(ii)}  \quad & \lim_{k \longrightarrow \infty} \left( \,perimeter(\tilde{\wp_k})\right) = \sqrt{A' \, (4 \pi -\kappa\, A')}.
\end{align*}
\noindent {\it Proof.}\\
\textbf{(i)}\quad 
\noindent \underline{$\kappa=0$} : \quad By (29), $r^2_k = r(k, A_k)^2=  \dfrac{2\,A_k}{k \, \sin \left( \frac{2\, \pi}{k}\right) }$. Therefore,
$$\lim_{k \longrightarrow \infty}r^2_k = \lim_{k \longrightarrow \infty}  \dfrac{2\,A_k}{2\,\pi}\dfrac{1}{\dfrac{\sin \left( \frac{2\, \pi}{k}\right)}{\left( \frac{2\, \pi}{k}\right) } }=\dfrac{A'}{\pi}.$$
Thus $$\lim_{k \longrightarrow \infty}r_k = \sqrt{\dfrac{A'}{\pi}}.$$
\noindent \underline{$\kappa \neq 0$} : \quad By (24), 
$$\dfrac{\cot\left( \frac{\pi}{k}\right)}{C_\kappa(r_k)}=\tan \left( \dfrac{A_k + \kappa \, (k-2)\,\pi }{2\,\kappa\,k}\right)=\tan \left( \dfrac{\kappa \,A_k- 2\pi}{2\,k} +\dfrac{\pi }{2}\right) = \cot \left( \dfrac{2\pi-\kappa \,A_k }{2\,k} \right).
$$
Hence, $$C_\kappa(r_k)=\cot\left( \dfrac{\pi}{k}\right)\,\tan \left( \dfrac{2\pi-\kappa \,A_k }{2\,k} \right).  
$$
Therefore, 
\begin{eqnarray*}
\kappa\, S_\kappa^2(r_k)
&=& 1- C_\kappa^2(r_k)
=\dfrac{\sin^2\left( \frac{\pi}{k}\right)\,\cos^2 \left( \frac{2\pi-\kappa \,A_k }{2\,k} \right)-\cos^2\left( \frac{\pi}{k}\right)\,\sin^2 \left( \frac{2\pi-\kappa \,A_k }{2\,k} \right) }{\sin^2\left( \frac{\pi}{k}\right)\,\cos^2 \left( \frac{2\pi-\kappa \,A_k }{2\,k} \right)}\\
&=&\dfrac{\sin\left( \frac{4\,\pi-\kappa\, A_k}{2\,k}\right)\,\sin \left( \frac{\kappa \,A_k }{2\,k} \right)}{\sin^2\left( \frac{\pi}{k}\right)\,\cos^2 \left( \frac{2\pi-\kappa \,A_k }{2\,k} \right)}.
\end{eqnarray*}
That is,
\begin{eqnarray*}
 S_\kappa^2(r_k)
&=& \dfrac{\kappa\,\sin \left( \frac{\kappa \,A_k }{2\,k} \right)\,\sin\left( \frac{4\,\pi-\kappa\, A_k}{2\,k}\right)}{\sin^2\left( \frac{\pi}{k}\right)\,\cos^2 \left( \frac{2\pi-\kappa \,A_k }{2\,k} \right)}\\
&=& \dfrac{\sin \left( \frac{A_k }{2\,k} \right)\,\sin\left( \frac{4\,\pi-\kappa\, A_k}{2\,k}\right)}{\sin^2\left( \frac{\pi}{k}\right)\,\cos^2 \left( \frac{2\pi-\kappa \,A_k }{2\,k} \right)}.
 \end{eqnarray*}
So,
$$\lim_{k\longrightarrow \infty}S_\kappa^2(r_k)=\dfrac{\left(4 \pi - \kappa\,A'\right)\,A'}{4 \pi^2}.$$
Hence, 
$$\lim_{k\longrightarrow \infty} r_k = AS_\kappa \left( \dfrac{ \sqrt{ \left( 4 \pi - \kappa\,A'\right)\,A'}}{2 \, \pi}\right).$$
\textbf{(ii)}\quad Put $r_0 = AS_\kappa \left( \sqrt{A' \, (4 \pi - \kappa\,A')}/(2\, \pi)\right)$. Now, each
$\tilde{\wp_k}$ is a regular $k$-gon inscribed in a circle of radius $r_k$ in $M_\kappa^2$, and $(r_k) \longrightarrow r_0$
as $k \longrightarrow \infty$ by \textbf{(i)}. Hence, $\left(perimeter(\tilde{\wp_k}) \right)$ converges to the perimeter
of the circle of radius $r_0$ in $M_\kappa^2$.
$$ ~~~~~~~~~~~~~~~~~~~\therefore ~ \lim_{k \longrightarrow \infty} \left(\, perimeter(\tilde{\wp_k}) \right)=
2 \pi\, S_\kappa( r_0 )= \sqrt{A' \, (4 \pi - \kappa A')}.\hfill ~~~~~~~~~~~~~~~~~~~~~~~~\blacksquare$$
\begin{center}
4. ISOPERIMETRIC PROBLEM FOR POLYGONS IN $M_\kappa^2$ 
\end{center} 
\noindent \textbf{Proof of Theorem 1 :}\quad Fix $n \geq 3$ in $\N ~ \& \; A \in \begin{cases} 
(0, \infty) & \text{if $\kappa=0$},\\
(0,2\pi) & \text{if $\kappa=1$},\\       
(0,(n-2)\,\pi) & \text{if $\kappa=-1$}.
\end{cases}$
 Let $\F$ be the family of all polygons with $n$ vertices in $M_\kappa^2$  having area at least $A$. By Proposition 3.2, there exists a regular $n$-gon $\wp_{n,r}$ of area equal to $A$. So, $\F$ is
a nonempty family. Define $ L = glb \left\{  perimeter(\wp) \, :\, \wp \in \F \right\}$. By Lemma 2.11, $perimeter(\wp_{n,r}) < 2 \pi $ if $\kappa =1$. Hence, $L < 2 \pi $ if $\kappa =1$. Let $\left( \wp_k \right)_{k \in \N}$ be
a sequence in $\F$ such that $\left( perimeter(\wp_k) \right) \searrow L$ as $k \longrightarrow \infty$ and $perimeter (\wp_k) < 2 \pi \; \forall \; k $ if $\kappa=1$. We assume $p_0$ as in (1) is a vertex of $\wp_k \; \forall \; k \in \N$. By Lemma 2.14, we can assume that if $\kappa =1$ then each $\wp_k$ is contained in the open hemisphere $B_1(p_0,\pi/2)$. Let $X_k^{(1)}, \,  X_k^{(2)}, \, \ldots, \, X_k^{(n)}$ be the vertices of $\wp_k$ with $X_k^{(1)}=p_0$ occurring in a cyclic order\,(determined by `boundary orientation' of $\partial \wp_k$) for all $k \in \N$. Then without loss of generality $\wp_k \in B_\kappa(p_0,L+1)\; \forall \; k$ when $\kappa \neq 1$. As $M_1^2$ is a compact manifold and $\overline{B_\kappa(p_0,L+1)}$ is compact in $M_\kappa^2$ for $\kappa \neq 1$, each sequence $\left( X_k^{(j)}\right)_{k \in \N}$ admits a converging subsequence $\forall \; j = 1, \ldots, n$.
Thus, without loss of generality, we can assume that $\left(
 X_k^{(j)}\right)_{k \in \N}$ converges to some
$Y_j$ in $M_\kappa^2
$, $ \forall \; j = 1, \ldots, n$. Clearly, $Y_1=p_0$. Let $Y_{n+1}:=Y_1$. Then $\cup_{i=1}^{n}[Y_i,Y_{i+1}]$ is a simple closed curve in $M_\kappa^2$. 

When $\kappa =1$, $L=\sum_{i=1}^{n}d_1(Y_i, Y_{i+1})< 2\, \pi$. and hence by Lemma 2.14, there exists a polygon $\wp_0$ contained in an open hemisphere 
of $M_1^2$ having $ Y_1,\, Y_2, \ldots, Y_n$ as its vertices occurring in a cyclic order. 
In particular, $\wp_0$ is a proper polygon when $\kappa=1$. When $\kappa \neq 1$, let $\wp_0$ be the polygon in $M_\kappa^2$ with $\cup_{i=1}^{n}[Y_i,Y_{i+1}]$ as its boundary. Then
$perimeter(\wp_0) = L$ and $area(\wp_0) \geq A$. It remains to show that $\wp_0$ is a regular $n$-gon.\\[1mm]
\noindent $\boldsymbol{\wp_0}$ {\it \textbf{is convex}} \textbf{:} \quad If not, by Theorem 2.9, $\exists \, j \in \{ 1, \ldots, n\}$ such that $\wp_0$
 is not locally convex at the vertex $Y_j$. Put $Y_0 : = Y_n$ and $Y_{n+1}:= Y_1$. Then we can choose a point
 $Y'_j$ on the side $[Y_{j-1},Y_j]$ such that the triangle $[Y'_j, Y_j, Y_{j+1}]$
 does not intersect interior of $\wp_0$. Now consider polygon $\wp$ in $M_\kappa^2$ having $Y_1,\, Y_2, \ldots,Y_{j-1}, \, Y'_j, \, Y_{j+1}, \ldots,\, Y_n$ as its vertices occurring in a cyclic order.
 Then $perimeter(\wp) < perimeter(\wp_0)$ and $area(\wp) > area(\wp_0) \geq A$. This contradicts the fact that
 $\wp_0$ is a perimeter minimizer in $\F$. So, $\wp_0$ is a convex $n$-gon.\\[1mm]
\noindent $\boldsymbol{area(\wp_0) = A} $ \textbf{:} \quad Suppose $area(\wp_0) = A + \delta$ with $\delta >0$. Choose a point $Y'_2 \in [Y_2, Y_3]$ such that $Y'_2 \not \in \{Y_2, Y_3\}$ and area of the triangle $[Y_1, Y_2, Y'_2]$
is less than $\delta$. As $\wp_0$ is convex the triangle $[Y_1, Y_2, Y'_2]$ is contained in $\wp_0$.
Then the polygon with vertices $ Y_1,\, Y'_2,\, Y_3, \ldots, Y_n$ occurring in a cyclic order has area greater than $A$ and perimeter
less than that of $\wp_0$. This is not possible. So, $area(\wp_0) = A $.\\[1mm]
\noindent $\boldsymbol{\wp_0}$ {\it \textbf{is equilateral}} \textbf{:} \quad If $\wp_0$ is not equilateral, then there exists
two successive sides of $\wp_0$ which are of unequal lengths. Suppose $d_\kappa(Y_1,Y_2)=:b \neq c:= d_\kappa(Y_2, Y_3)$.
Join $Y_1$ and $Y_3$ by the geodesic segment $[Y_1,Y_3]$. 
The triangle
$[Y_1,Y_2,Y_3]$ is contained in $\wp_0$. The geodesic segment $[Y_1,Y_3]$ as above divides $\wp_0$ in two polygons, namely, triangle $[Y_1,Y_2,Y_3]$ and the $(n-2)$-gon with vertices $Y_1, Y_3, Y_4, \ldots, Y_n $ occurring in a cyclic order. We call this $(n-2)$-gon as $\wp$. By Proposition 2.8, $\exists$ an isosceles triangle $[Y_1,Y'_2,Y_3]$ in $M_\kappa^2$
such that $Y'_2, \, Y_2$ lie on the same half-space whose boundary contains $[Y_1,Y_3]$, $ d_\kappa(Y_1, Y'_2) = (b+c)/2 = d_\kappa(Y'_2, Y_3)$ and $area([Y_1, Y'_2, Y_3]) > area([Y_1, Y_2, Y_3])$. Then the polygon $\wp \cup [Y_1,Y'_2,Y_3]$ 
is a perimeter minimizer in $\F$ with
area strictly greater than $A$. This is not possible as seen in the previous step. Thus $\wp_0$ is an equilateral $n$-gon.\\[1mm]
Let `$a$' denote the side length of the convex equilateral $n$-gon $\wp_0$.\\
\noindent $\boldsymbol{\wp_0}$ {\it \textbf{is equiangular}} \textbf{:} \quad We prove this by considering $n$ even, $n$ odd cases separately.\\[1mm]
\noindent $\boldsymbol{n}$ \textbf{is even :} \quad Let $n = 2\,k$. By Lemma 2.11, $k a = \frac{n}{2}\, a < \pi$ if $\kappa=1$. Join $Y_1$ to $Y_{1+k}$
by the geodesic segment $[Y_1,Y_{1+k}]$ contained in $\wp_0$. Let $O$ be the mid-point of $[Y_1,Y_{1+k}]$.
Enough to show that $d_\kappa(O, Y_i)= r :=d_\kappa(O,Y_1) \; \forall \; i = 2, \ldots,n$.
\par The segment $[Y_1,Y_{1+k}]$ divides $\wp_0$ into two convex polygons $\wp_1$, $\wp_2$ with $k+1$ vertices.
Let $\rho$ be the reflection of $M_\kappa^2$ through the {\it line} containing $[Y_1,Y_{1+k}]$.
If $area(\wp_1) > area(\wp_2)$, then $\wp_1 \cup \rho(\wp_1)$ gives a polygon of perimeter
$L$ and area greater than $A$. This is not possible. So, $area(\wp_1) = area(\wp_2)$.
\par Consider the family of all the convex  polygons with $(k+1)$ vertices in $M_\kappa^2$ whose all sides but one are of equal length $a$. Let $\wp_1' = [Y'_1, \ldots,Y'_{k+1}]$ be the area maximizer in this family with
$[Y'_1, Y'_{k+1}]$ being the `remaining side'. Then $area(\wp_1') \geq area(\wp_1)$.
If $area(\wp_1') > area(\wp_1)$, then we can produce a polygon of perimeter $L = 2 k a$ and
area greater than $A$ by reflecting $\wp_1'$ through the side $[Y'_1, Y'_{k+1}]$. Hence
$area(\wp_1') = area(\wp_1)$. Now, by lemma 2.13 it follows that
 $d_\kappa(O, Y_i)= r \; \forall \; i = 1, \ldots,k+1$. Similarly one can prove that
 $d_\kappa(O, Y_i)= r \; \forall \; i = k+1, \ldots,2k$.\\[1mm]
\noindent $\boldsymbol{n}$ \textbf{is odd :} \quad Suppose $\wp_0$ is not equiangular. Let $\wp_R^n$ 
be a regular $n$-gon of side $a$ in $M_\kappa^2$. By Proposition 2.6, there exists an isosceles triangle $T$ in $M_\kappa^2$ having base $a$
and very small angles $\alpha$ at the base. Let $\wp_0^{2n }$ be the polygon with $2n$ sides obtained by `pasting' triangle congruent to $T$ on each side of $\wp_0$ so that $ area(\wp_0^ {2n}) = area(\wp_0) + n \; area(T)$. This is possible since $\wp_0$ is a convex polygon. As $\wp_0$ is not equiangular, $\wp_0^ {2n}$ is not equiangular.
 Similarly construct a regular $2n$-gon $\wp_R^{2n}$ by `pasting' triangle congruent to $T$ on each side of $\wp_R^n$.
 By the `$n$-even' case, as perimeters of $\wp_0^{2n}$ and $\wp_R^{2n}$ are equal, $area(\wp_0^ {2n})   < area(\wp_R^{2n})$. Therefore, $area(\wp_0) + n \; area(T)  < area(\wp_R^n) + n \; area(T)$. So, $A =   area(\wp_0)   < area(\wp_R^n)$.
Also, $\, perimeter(\wp_R^n) = n\,a = perimeter(\wp_0)$.
Thus $\wp_R^n \in \F$  and $\wp_R^n$ is a `perimeter minimizer' in $\F$. We have a contradiction as any
`perimeter minimizer' in $\F$ has area $A$. We conclude that $\wp_0$ is equiangular. \hfill $\blacksquare$
\begin{center}
5. THE ISOPERIMETRIC PROBLEM IN $M_\kappa^2$ 
\end{center}
\noindent \textbf{Notations :}\quad For a piecewise smooth simple closed curve $\gamma$ in $M_\kappa^2$ let $\ell(\gamma)$ denote the arc-length of $\gamma$. For $\kappa =1$ if such a curve $\gamma$ lies in a hemisphere $S^+$ then $\gamma$ encloses a domain $D_\gamma$ contained in $S^+$. If $\kappa \neq 1$ then such a curve $\gamma$ always encloses a unique relatively compact domain $D_\gamma$ contained in $M_\kappa^2$. We denote $area(D_\gamma)$ by $A(\gamma)$.\\[1mm]
\noindent \textbf{Proof of Theorem 2 :}\\
\noindent \textbf{Case (i)} \quad $\kappa=1$ and $A= 2 \pi$ :\quad {\it There exists a unique perimeter minimizer among all piecewise smooth simple closed curves in} $M_1^2$ {\it enclosing area} $2\pi$, {\it and it is a great circle }:\\[1mm]
Let $\J$ be the family of all piecewise smooth simple closed curves in $M_1^2$ enclosing area $2\pi$. Let $S^2_+ : = \{(x,y,z)\in S^2 \;|\; z \geq 0 \}$. Since $\partial S^2_+ \in \J$, $\J \neq \emptyset$. 
\par If $\exists \;\C \in \J$ with $\ell(\C)<2\pi$ then by Lemma 2.14, $\C$ is contained in an open hemisphere. This contradicts the fact that $\C$ encloses area $2\pi$
. Hence, $$\ell(\C)\geq 2\pi \; \forall \; \C \in \J.\eqno(30)$$ 
$$\mbox{Define}~ L:= \inf \{ \ell(\C)\;|\; \C \in \J \}.$$
Since $\partial S^2_+ \in \J$, $$L\leq \ell(\partial S^2_+)=2\pi.\eqno(31)$$ Thus from $(30)$ and $(31)$, we get $L=2 \pi$ and $\partial S^2_+$ is a perimeter minimizer over $\J$. 

Let $\C_0$ be a perimeter minimizer over $\J$. That is, $\ell(\C_0)=L =2 \pi$ and $A(\C_0) := \mbox{area enclosed by } \C_0=2\pi$. If $\C_0$ is not boundary of a digon 
then by Lemma 2.15, $\C_0$ is contained in an open hemisphere, a contradiction again. Hence $\C_0 = \partial D_{x,\alpha}$ for some $x \in M_1^2$ and $\alpha \in [0, \pi]$. Therefore, $2 \pi = A(\C_0) = area(D_{x,\alpha}) = 2 \alpha$. This implies that $\alpha =\pi$. Thus $\C_0 = \partial D_{x,\pi}$, that is a great circle.\\[1mm] 
\noindent \textbf{case (ii)}\quad $\kappa \in \{-1,0,1\}$ and $A <2 \pi$ if $\kappa=1$ :\\
Let $p_0 \in M_\kappa^2$ be as in (1). For $r_0>0$ ($r_0 <\pi/2$ if $\kappa=1$), the circle $\mathcal{C}_{r_0}:=\partial B_\kappa(p_0,r_0)$ encloses a domain of area $4 \pi S_\kappa^2\left( \frac{r_0}{2}\right) $ ($< 2 \pi$ if $\kappa=1$) with perimeter $2 \pi S_\kappa( r_0)$ ($< 2 \pi$ if $\kappa=1$). So, for $\kappa=1$ we need to consider piecewise smooth simple closed curves of lengths strictly less than $2 \pi$ only. By Lemma 2.14, any such curve lies in a hemisphere.\\[1.75mm]
\noindent \textbf{Step 1.} (Existence)\quad {\it Among all piecewise smooth simple closed curves in}
 $M_\kappa^2$ {\it enclosing area} $A$, {\it a circle of radius}
 $AS_\kappa \left(\sqrt{A \, (4 \pi -\kappa\, A)}/(2 \pi)\right)$ {\it in} $M_\kappa^2$ {\it has least perimeter }:\\[1mm]
Let $\mathcal{G}$ denote the family of all piecewise smooth simple closed curves in $M_\kappa^2$ (in $S^2_+ := \left\{  (x, y, z) \in S^2 \, | \, z \geq 0 \right\}$ if $\kappa=1$) enclosing area at least $A$. Let $\mathcal{C} \in \mathcal{G}$ be arbitrary. If $Y_1,\, Y_2, \ldots, Y_n$ are points on a curve $\C \in \G$ which appear in a cyclic order [with
$d_\kappa(Y_i,Y_{i+1}) < \pi \; \forall \; i = 1, \ldots , n$ if $\kappa =1$ ($Y_{n+1}:= Y_1$)], then they determine a polygon with vertices
$Y_1,\, Y_2, \ldots, Y_n$ (which is contained in $S^2_+$ if $\kappa=1$). Define $L:= glb \left\{  \ell(\C) \; | \; \C \in \G \right\}$. Let $\left(\C_n\right)_{n \in \N}$ be a sequence in $\G$ such that $\ell(\C_n) \searrow L$
as $n \longrightarrow \infty$. Let $p_0 \in M_\kappa^2$ be as in (1). We may assume that $ p_0 \in \C_{n} $ and that  $ l(\C_{n}) \leq L+1 ,~\forall\ n \in \mathbb{N}$. Hence $ \C_{n} \subset B_\kappa(p_0,L+1),~\forall\ n \in \mathbb{N} $. Then $ A(\C_{n}) \leq A\left(B_\kappa(p_0,L+1)\right)=4\pi S_\kappa^2\left( \frac{L+1}{2}\right) ~ \forall\ n \in \mathbb{N} $. Therefore, we may assume, after taking
a subsequence of $\left(\C_n\right)_{n \in \N}$ if necessary, that
$$ \ell(\C_n) \searrow L \; \mbox{ and} \; \left( A(\C_n)\right)_{n \in \N} \longrightarrow  : A' \geq A. \eqno{(32)}$$
As each $\C_n$ is a piecewise smooth simple closed curve, we can approximate $\C_n$ by the boundary of a polygon $\wp_{k(n)}$ in $M_\kappa^2$
with $k(n)$ sides : i.e., vertices of $\wp_{k(n)}$ lie on $\C_n$, $0 < \ell(\C_n) - \ell(\partial \wp_{k(n)}) < \frac{1}{n}$ and $|area(\wp_{k(n)}) - A(\C_n)| < \frac{1}{n} \;\; \forall \; n \in \N\,$. (Here $\partial \wp$ denotes the boundary of $\wp$). Then it follows that 
$$ \lim_{n \longrightarrow \infty} \ell(\partial \wp_{k(n)}) = L ~~\mbox{ and }\lim_{n \longrightarrow \infty} area(\wp_{k(n)}) = A'.$$ 
Put $A_n := area(\wp_{k(n)})  \; \forall \; n \in \N$. By Proposition 3.2, for each $n \in \N$ there exists regular
$k(n)$-gon $\tilde{\wp}_{k(n)}$ in $M_\kappa^2$ of area $A_n$. By Theorem 1,
$\ell(\partial \tilde{\wp}_{k(n)}) \leq \ell(\wp_{k(n)}) \; \forall \; n \in \N$. Let $r_{k(n)}$ be the radius
of the circle in $M_\kappa^2$ in which $\tilde{\wp}_{k(n)}$ is inscribed. By Corollary 3.3, 
$$\lim_{n \longrightarrow \infty} r_{k(n)} = AS_\kappa \left(\sqrt{A' \, (4 \pi - \kappa\,A')}/(2 \pi)\right) =: r_0.$$ 
Thus, 
$$A' = 4 \pi S_\kappa^2 \left( \dfrac{r_0}{2}\right). \eqno{(33)}$$
Again by Corollary 3.3, $\lim_{n \longrightarrow \infty} \ell(\partial \tilde{\wp}_{k(n)}) = \sqrt{A' \, (4 \pi -\kappa\, A')} = 2 \pi S_\kappa(r_0)$. Note that $\ell(\partial \tilde{\wp}_{k(n)})$ $\leq \ell(\partial \wp_{k(n)}) \leq \ell(\C_n) \; \forall \; n \in \N$. Therefore,
$$  L = \lim_{n \longrightarrow \infty}\ell(\C_n) \geq
\lim_{n \longrightarrow \infty} \ell(\partial \tilde{\wp}_{k(n)}) = 2 \pi S_\kappa(r_0). \eqno{(34)}$$
Let $\C_{r_0}$ denote the circle in $M_\kappa^2$ (in $S^2_+$ if $\kappa=1$) of radius $r_0$. Then $A(\C_{r_0}) = 4 \pi S_\kappa^2\left(  \frac{r_0}{2}\right) = A'$ (by (33)). So,
by (32), $\C_{r_0} \in \G$, and by the definition of $L$,
$$ L \leq \ell(\C_{r_0}) = 2 \pi S_\kappa( r_0). \eqno{(35)}$$
By $(34) ~ \& ~ (35)$, $L =  2 \pi S_\kappa( r_0) = \ell(\C_{r_0})$. Hence $\C_{r_0}$ is a perimeter minimizer in $\G$.\\[1mm]
Finally we show that $A(\C_{r_0}) = A$ : \quad If $A(\C_{r_0}) \neq A$ then by $(32) ~ \& ~ (33)$,
$A(\C_{r_0}) = 4 \pi S_\kappa^2\left(  \frac{r_0}{2}\right)= A' > A$. Then we can replace a small portion of circle $\C_{r_0}$ by a geodesic arc
and produce a curve $\tilde{\C}$ in $M_\kappa^2$ (in $S^2_+$ if $\kappa=1$) with $\ell(\tilde{\C}) < \ell(\C_{r_0}) = L$ and $A < A(\tilde{C}) < A(\C_{r_0})$. Then $\tilde{\C} \in \G$ with $\ell(\tilde{\C}) < L$. This is not possible. Hence $A(\C_{r_0}) = A$.
 
\noindent \textbf{Step 2.} (Uniqueness)\quad {\it Among all piecewise smooth simple closed curves in}
 $M_\kappa^2$ {\it enclosing area} $A$, {\it any perimeter minimizer is a circle in\;} $M_\kappa^2$ {\it of radius}\\
 $AS_\kappa \left(\sqrt{A \, (4 \pi - \kappa\,A)}/(2 \pi)\right)$ :\\[1mm]
Consider the family $\G$ of all piecewise smooth simple closed curves in $M_\kappa^2$ (in $S^2_+$ if $\kappa=1$) enclosing area at least $A$ and of perimeter strictly less than $2 \pi$ for $\kappa=1$. Put
$L:= glb \left\{  \ell(\C) \; | \; \C \in \G \right\}$. In Step 1 above, we have proved the existence of curve in $\G$
which is a perimeter minimizer. Let $\C_0 \in \G $ be any perimeter minimizer. Then $\ell(\C_0) =L$ ($< 2 \pi$ if $\kappa=1$).
Let $D_0$ be the domain in $M_\kappa^2$ (in $S^2_+$ if $\kappa=1$) enclosed by $\C_0$. By the arguments similar to those made in the proof of
 Theorem 1, we can show that $D_0$ is convex and $area(D_0) = A$.
\par Fix a point $P$ on $\C_0$. Let $Q$ be the point on $\C_0$ which divides $\C_0$ into two arcs $\C_0^+$, $\C_0^-$
of equal length. As $\ell(\C_0) < 2 \pi$ in $M_1^2$ , $Q \neq -P$ if $\kappa=1$. Let $[P,Q]$ denote the geodesic segment joining $P \; \& \; Q$ in
$D_0$. This segment divides $D_0$ into two regions $D_0^+ \; \& \; D_0^-$. If $area(D_0^+) < area(D_0^-)$, then consider
$\tilde{D}_0 := D_0^- \cup \rho(D_0^-)$ where $\rho$ is the reflection in $M_\kappa^2$ through the {\it line} containing $[P,Q]$.
Then boundary $\tilde{\C}_0$ of $\tilde{D}_0$ is a perimeter minimizer in $\G$ and $area(\tilde{D}_0) > area(D_0)$.
This is not possible. Hence, $[P, Q]$ divides $D_0$ into two regions of equal area.
\par Let $O$ be the mid-point of $[P,Q]$ and $r_0 := d_\kappa(P,Q)/2$. We show that
$d_\kappa(O,M) = r_0 \; \forall \; M \in \C_0$ : Suppose $\exists \, M \in \C_0$ such that $d_\kappa(O,M) \neq r_0$.
 Let $D_0^+$ be the region containing $M$ with $\C_0^+ \cup [P,Q]$ as its boundary.
As $D_0$ is convex, the triangle $[P,M,Q]$ of $M_\kappa^2$ is contained in $D_0^+$. Now, $d_\kappa(P,M) + d_\kappa(M,Q) \leq \ell(\C_0^+) = L/2$ ($ < \pi$ if $\kappa=1$) and $M$ does not lie on the circle in $M_\kappa^2$ of radius $r_0$ and center $O$.
\par By the arguments similar to those made in the proof of Lemma 2.13, we can construct a domain $\widetilde{D_0^+}$ in $M_\kappa^2$ (in $S^2_+$ if $\kappa=1$) of area strictly bigger than $area(D_0^+)$ whose boundary consists of a curve\;
$\widetilde{\C_0^+}$ which is congruent to \;$\C_0^+$ and a geodesic segment $[P',Q']$\,($P',\, Q'$ are the endpoints of \;$\widetilde{\C_0^+}$). Reflecting $\widetilde{D_0^+}$ through the {\it line} containing $[P', Q']$ we can produce a domain
 $\widetilde{D_0}$ of area strictly bigger than $A$ and perimeter of boundary of $\widetilde{D_0}$ equal to $L$. This is
 not possible. Hence $d_\kappa(O, M) = r_0$ and $\C_0 = \partial B_\kappa(O, r_0)$. \hfill $\blacksquare$\\[3mm]
\noindent \textbf{Proof of Corollary 3} \quad Let $\C$ be a piecewise smooth simple closed curve  having $m$ components each enclosing area $A_i>0$. Let $r_i :=AS_\kappa \left(\sqrt{A_i \, (4 \pi - \kappa\, A_i)}/(2 \pi)\right)$ ($1 \leq i \leq m$). Let $\tilde{\C}$ denote the disjoint union of the circles $\tilde{\C}_i$ of radius $r_i$, $1 \leq i \leq m$. Applying Theorem 2 to each component of $\C$ we get that $perimeter(\tilde{\C})\leq perimeter(\C)$. Now, it is easy to see that a single circle with radius $AS_\kappa \left(\frac{\sqrt{A \, (4 \pi - \kappa\, A)}}{2 \pi}\right)$ is the best. \hfill $\blacksquare$\\[3mm]   
\noindent \textbf{Remarks :}
\begin{enumerate}
\item Fix $L_0 \in (0, 2\pi]$. Put $r_0:= \arcsin \left(L_0/(2\pi)\right)\in (0, 2\pi]$ and $A_0:= 4 \pi \sin^2 \left( \frac{r_0}{2}\right)$. Let $\C$ be any piecewise smooth simple closed curve in $M_1^2$ having arc-length $\ell(\C)= L_0$. From Theorem 2, it follows that among all such curves, area maximizer is the circle $\C_{r_0}$. For, consider the family $\mathcal{F} = \{\,$all piecewise smooth simple closed curves in $M_1^2$ enclosing area $\geq A_0 \}$. If $A(\C)\geq A_0$, then $\C \in \mathcal{F}$ and $\ell(\C)=L_0= \ell(\C_{r_0})= \inf \,\{\, \ell(\tilde{\C})\,|\, \tilde{\C}\in \mathcal{F}\,\}$. By Theorem 2, $\C = \C_{r_0}$ and $A(\C)= A_0$. 
\item A shorter though less elementary approach to prove Theorem 1 for $M_{-1}^2$ is to first prove Theorem 2 for this case and then derive the results for $n$-gons using Heron's formula or L'Huilier's Theorem as in [\cite{Quinn}, Proposition 2.15]. 
\end{enumerate}

\noindent \textbf{Proof of Theorem 4 :}\quad Let $\C$ be any piecewise smooth simple closed curve in $M_\kappa^2$ with arc-length $\ell:=\ell(\C)$ and enclosing area $A:= A(\C) >0$ ($A \leq 2\pi$ if $\kappa=1$).\\
\quad \textbf{Case (i)} \quad $\kappa=1$ and $A=2 \pi$ :\\
 Let $\J$ and $L$ be as in the proof of Theorem 2 for the corresponding case. 
 Also recall that $L=2\pi$. Therefore, $L^2 = 4 \pi ^2 = 4 \pi A - A^2$ and hence $\ell^{\, 2} =[\ell(\C)]^2 \geq L^2 = 4 \pi A - A^2$  holds for all $\C \in \J$.
\par If for a curve $\C$ in $\J$, $[\ell(\C)]^{2} = 4 \pi A - A^2 = 4 \pi ^2$, then $\ell(\C) = \sqrt{4 \pi A - A^2}=2 \pi = L$ and $\C$ is a perimeter minimizer in $\J$.
 By Theorem 2, $\C$ is a great circle i.e., a circle in $M_1^2$ of radius $ \frac{\pi}{2} = \arcsin \left(\sqrt{A \, (4 \pi - A)}/(2 \pi)\right) $.
\\[1mm] 
\noindent \textbf{Case (ii)} \quad $\kappa \in \{-1,0,1\}$ and $A <2\pi$ if $\kappa=1$:\\
 Let $\G$ and $L$ be as in the proof of Theorem 2 for the corresponding case. 
 Put
$r_0 := AS_\kappa \left(\sqrt{A \, (4 \pi -\kappa\, A)}/(2 \pi)\right)$. By Theorem 2,
the circle $\C_{r_0}$ of radius $r_0$ is the unique perimeter minimizer in $\G$. Therefore, $L^2 = \left(2 \pi S_\kappa( r_0) \right)^2 = 4 \pi A - \kappa\, A^2$
and hence $\ell^{\, 2} =[\ell(\C)]^2 \geq L^2 = 4 \pi A - \kappa\,A^2$  holds for all $\C$ in $\G$.
\par If for a curve $\C$ in $\G$, $[\ell(\C)]^2= 4 \pi A - \kappa\,A^2$, then $\ell(\C) = \sqrt{4 \pi A - \kappa\,A^2} = L$ and $\C$ is a perimeter minimizer in $\G$.
 By Theorem 2, $\C$ is a circle of radius $r_0$ in $M_\kappa^2$. \hfill $\blacksquare$

\begin{center}
6. APPENDIX
\end{center}
\noindent \textbf{\large 6.1 Appendix A}\\
We state some formulae about $S_\kappa$ and $C_\kappa$ when $\kappa \neq 0$.
$$C_\kappa (-a)= C_\kappa (a), ~~S_\kappa(-a)= -S_\kappa(a).~~~~~~~~~~~~~~~~~~~~~~~~~~~~~~~~~~~~~~~~~~~~~~~~~~~~ \eqno{(A-1)}$$
$$S_\kappa (a+b)= S_\kappa (a)\, C_\kappa (b) + C_\kappa (a)\, S_\kappa (b).~~~~~~~~~~~~~~~~~~~~~~~~~~~~~~~~~~~~~~~~~~~~~~~~~~~~\eqno{(A-2)}$$
$$S_\kappa (a-b)= S_\kappa (a)\, C_\kappa (b) - C_\kappa (a)\, S_\kappa (b).~~~~~~~~~~~~~~~~~~~~~~~~~~~~~~~~~~~~~~~~~~~~~~~~~~~~\eqno{(A-3)}$$
$$S_\kappa (2\, a)= 2\, S_\kappa (a)\, C_\kappa (a).~~~~~~~~~~~~~~~~~~~~~~~~~~~~~~~~~~~~~~~~~~~~~~~~~~~~~~~~~~~~~~~~~~\eqno{(A-4)}$$
$$C_\kappa^2(a) = 1- \kappa \, S_\kappa^2(a).~~~~~~~~~~~~~~~~~~~~~~~~~~~~~~~~~~~~~~~~~~~~~~~~~~~~~~~~~~~~~~~~~~~~~~~~~~\eqno{(A-5)}$$
$$C_\kappa (a+b)= C_\kappa (a)\, C_\kappa (b) -\kappa\, S_\kappa (a)\, S_\kappa (b).~~~~~~~~~~~~~~~~~~~~~~~~~~~~~~~~~~~~~~~~~~~~~~~~~~~~~~~~~~~~~~~~~~~\eqno{(A-6)}$$
$$C_\kappa (a-b)= C_\kappa (a)\, C_\kappa (b) +\kappa\, S_\kappa (a)\, S_\kappa (b).~~~~~~~~~~~~~~~~~~~~~~~~~~~~~~~~~~~~~~~~~~~~~~~~~~~~~~~~~~~~~~~~~~~\eqno{(A-7)}$$
$$C_\kappa (2\,a)= C_\kappa^2 (a) -\kappa\, S_\kappa^2 (a)=1-2\,\kappa\, S_\kappa^2 (a)= 2\,C_\kappa^2 (a) -1.~~~~~~~~~~~~~~~~~~~~~~~~~~~~~~~~~\eqno{(A-8)}$$
$$1-C_\kappa(a)= 2\, \kappa \, S_\kappa^2 (a/2).~~~~~~~~~~~~~~~~~~~~~~~~~~~~~~~~~~~~~~~~~~~~~~~~~~~~~~~~~~~~~~~~~~~~~~~~~~~~~~~~~~~~\eqno{(A-9)}$$
$$1+C_\kappa(a)= 2 \, C_\kappa^2 (a/2).~~~~~~~~~~~~~~~~~~~~~~~~~~~~~~~~~~~~~~~~~~~~~~~~~~~~~~~~~~~~~~~~~~~~~~~~~~~~~~~~~~~~\eqno{(A-10)}$$
$$C_\kappa(a+b)+ C_\kappa(a-b)= 2\, C_\kappa(a)\, C_\kappa(b).~~~~~~~~~~~~~~~~~~~~~~~~~~~~~~~~~~~~~~~~~~~~~~~~~~~~~~~~~~~~~~~~~~~\eqno{(A-11)} $$
$$C_\kappa(a+b)- C_\kappa(a-b)= -2\,\kappa\, S_\kappa(a)\, S_\kappa(b) ~~~~~~~~~~~~~~~~~~~~~~~~~~~~~~~~~~~~~~~~~~~~~~~~~~~~~~~~~~\eqno{(A-12)}$$
$$2 \, S_\kappa \left({\dfrac{a+b}{2}}\right)C_\kappa \left({\dfrac{a-b}{2}}\right) = S_\kappa (a) + S_\kappa (b).~~~~~~~~~~~~~~~~~~~~~~~~~~~~~~~~~~~~~~~~~~~~~~~~~~~~~\eqno{(A-13)} $$ 
$$ 2\,C_\kappa\left({\dfrac{a+b}{2}}\right)S_\kappa\left({\dfrac{a-b}{2}}\right) = S_\kappa(a) - S_\kappa (b).~~~~~~~~~~~~~~~~~~~~~~~~~~~~~~~~~~~~~~~~~~~~~~~~~~~~~ \eqno{(A-14)}$$
$$ 2\, C_\kappa\left({\dfrac{a+b}{2}}\right)C_\kappa\left({\dfrac{a-b}{2}}\right) = C_\kappa (a) + C_\kappa (b).~~~~~~~~~~~~~~~~~~~~~~~~~~~~~~~~~~~~~~~~~~~~~~~~~~~~~~~\eqno{(A-15)}$$
$$ -2\,\kappa\,S_\kappa\left({\dfrac{a+b}{2}}\right)S_\kappa\left({\dfrac{a-b}{2}}\right) = C_\kappa (a) - C_\kappa (b).~~~~~~~~~~~~~~~~~~~~~~~~~~~~~~~~~~~~~~~~~~~~\eqno{(A-16)}$$
\noindent \textbf{\large 6.2 Appendix B }\\[1mm]
\textbf{Trigonometric formulae for a triangle in} $\boldsymbol{M_\kappa^2}$ 
\\[1mm]
Let $[P,Q,R]$ be a triangle in  $ M_\kappa^{2} $ having angles $ \alpha, \beta, \gamma $ at its vertices and let $ a, b, c $ be sides opposite to angles $\alpha,\beta, \gamma$, respectively.
Put $ s = \dfrac{a+b+c}{2} .$ Then we have following formulae :\\
$$ \sin\dfrac{\gamma}{2}\, = \sqrt{\dfrac{S_\kappa(s-a) S_\kappa(s-b)}{S_\kappa (a) S_\kappa (b)}}.~~~~~~~~~~~~~~~~~~~~~~~~~~~~~~~~~~~~~~~~~~~~~~~~~~~~~~~~~~\eqno(B-1)$$
$$ \cos\dfrac{\gamma}{2}\, = \sqrt{\dfrac{S_\kappa (s)\,S_\kappa(s-c)}{S_\kappa (a)\, S_\kappa (b)}}.~~~~~~~~~~~~~~~~~~~~~~~~~~~~~~~~~~~~~~~~~~~~~~~~~~~~~~~~~~~~\eqno(B-2)$$
\textbf{The sine rule :}
$$\dfrac{\sin \alpha}{S_\kappa (a)} = \dfrac{\sin \beta}{S_\kappa (b)} = \dfrac{\sin \gamma}{S_\kappa (c)} =\dfrac{2 \sqrt{S_\kappa (s)\, S_\kappa(s-a) S_\kappa(s-b) S_\kappa(s-c)}}{S_\kappa (a)\, S_\kappa (b)\, S_\kappa (c)}.~~~~~~~~~~~~~\eqno(B-3)$$\\[1mm]
\noindent The following holds when $\kappa \neq 0$ :\\
$$ \sin\left({\dfrac{\alpha+\beta}{2}}\right) = \cos\dfrac{\gamma}{2}\,\, \dfrac{C_\kappa\left({\dfrac{a-b}{2}}\right)}{C_\kappa\left(\dfrac{c}{2}\right)\,}.~~~~~~~~~~~~~~~~~~~~~~~~~~~~~~~~~~~~~~~~~~~~~
~~~~~~~~~~~~~\eqno(B-4)$$ 
$$ \sin\left({\dfrac{\alpha-\beta}{2}}\right) = \cos\dfrac{\gamma}{2}\,\, \dfrac{S_\kappa\left({\dfrac{a-b}{2}}\right)}{S_\kappa\left( \dfrac{c}{2}\right) \,}.~~~~~~~~~~~~~~~~~~~~~~~~~~~~~~~~~~~~~~~~~~~~~~~
~~~~~~~~~~\eqno(B-5)$$
$$ \cos\left({\dfrac{\alpha+\beta}{2}}\right) = \sin\dfrac{\gamma}{2}\,\, \dfrac{C_\kappa\left({\dfrac{a+b}{2}}\right)}{C_\kappa\left(\dfrac{c}{2}\right)\,}.~~~~~~~~~~~~~~~~~~~~~~~~~~~~~~~~~~~~~~~~~~~~
~~~~~~~~~~~~~~~\eqno(B-6)$$ $$\cos\left({\dfrac{\alpha-\beta}{2}}\right) = \sin\dfrac{\gamma}{2}\,\,
\dfrac{S_\kappa\left({\dfrac{a+b}{2}}\right)}{S_\kappa \left(\dfrac{c}{2}\right)\,}.~~~~~~~~~~~~~~~~~~~~~~~~~~~~~~~~~~~~~~~~~~~~~~~~~~~\eqno(B-7)$$
\textbf{Proof of (B-1) :} \quad We give the proof of (B-1) for triangles in $M_\kappa^2$ ($\kappa \neq 0$). The proof for triangles in $M_0^2$ is similar and simpler.

By the Law of Cosine for triangles in $ M_\kappa^{2} $ ($\kappa \neq 0$) we have
$$\cos(\gamma) =  \dfrac{C_\kappa (c)-C_\kappa (a)\, C_\kappa (b)  }{\kappa\,S_\kappa (a)\,S_\kappa (b)}.\eqno{(B-8)}$$
By (A-9),
\begin{align*}
2 \sin^{2}\dfrac{\gamma}{2}\, &= 1 - \dfrac{C_\kappa (c)-C_\kappa (a)\, C_\kappa (b)  }{\kappa\,S_\kappa (a)\,S_\kappa (b)} =  \dfrac{\kappa\,S_\kappa (a) \,S_\kappa (b) + C_\kappa (a)\, C_\kappa (b) - C_\kappa (c)}{\kappa\,S_\kappa (a)\, S_\kappa (b)}
\\[1mm]
& =  \dfrac{C_\kappa(a-b)-C_\kappa (c)}{\kappa \,S_\kappa (a) \,S_\kappa (b)}~~~~~~~~~~~~~~~~~~~~~~~~~~~~~~~~~~~~~~~~~~~~~~~~~~~~~~~~\mbox{[by (A-7)]}\\[1mm]
&= \dfrac{2 \kappa \, S_\kappa\left({\dfrac{c+a-b}{2}}\right)S_\kappa\left({\dfrac{c-a+b}{2}}\right)}{\kappa\,S_\kappa (a)\,S_\kappa (b)}~~~~~~~~~~~~~~~~~~~[\mbox{by  (A-16) and (A-1)}]\\[1mm] &=\dfrac{2S_\kappa(s-b)S_\kappa(s-a)}{S_\kappa (a)\,S_\kappa (b)}.\end{align*}
Now (B-1) follows easily.\\[3mm]
\textbf{Proof of (B-2) :} \quad We give the proof of (B-2) for triangles in $M_\kappa^2$ ($\kappa\neq 0$). The proof for triangles in $M_0^2$ is similar and simpler.

By (B-8) and (A-10) we have :
\begin{eqnarray*} 2\cos^{2}\dfrac{\gamma}{2}\, &=& 1 + \dfrac{C_\kappa (c)-C_\kappa (a)\, C_\kappa (b)  }{\kappa\,S_\kappa (a)\,S_\kappa (b)}
=\dfrac {-C_\kappa (a)\,C_\kappa (b) + \kappa\,S_\kappa (a)\, S_\kappa (b) +C_\kappa (c)}{\kappa\, S_\kappa (a) \,S_\kappa (b}\\ 
&=& \dfrac{-C_\kappa(a+b) + C_\kappa (c) }{\kappa\,S_\kappa (a) \,S_\kappa (b)}~~~~~~~~~~~~~~~~~~~~~~~~~~~~~~~~~~~~~~~~~~~~~~~~~\mbox{[by (A-6)]}\\ &=& \dfrac{2 \kappa\,S_\kappa\left({\dfrac{a+b+c}{2}}\right) S_\kappa\left({\dfrac{a+b-c}{2}}\right)}{\kappa\,S_\kappa (a)\, S_\kappa (b)}~~~~~~~~~~~~~~~[\mbox{by  (A-16) and (A-1)}]
\\ &=& \dfrac{2 S_\kappa (s)\, S_\kappa(s-c)}{S_\kappa (a)\,S_\kappa (b)}.\end{eqnarray*} 
Now (B-2) follows easily.\\[3mm]
\textbf{Proof of (B-3) :} \quad By (A-4), (B-1) and (B-2) we get,
$$\dfrac{\sin\gamma}{S_\kappa (c)} = \dfrac{2\sin\dfrac{\gamma}{2}\, \cos\dfrac{\gamma}{2}}{S_\kappa (c)}=\dfrac{2\ \sqrt{S_\kappa (s)\, S_\kappa(s-a) S_\kappa(s-b) S_\kappa(s-c)}}{S_\kappa (a)\, S_\kappa (b)\, S_\kappa (c)}.~~~~~~~~~~~~$$ Hence,
$$\dfrac{\sin\alpha}{S_\kappa (a)} = \dfrac{\sin\beta}{S_\kappa (b)} = \dfrac{\sin\gamma}{S_\kappa (c)} = \dfrac{2 \sqrt{S_\kappa (s)\, S_\kappa(s-a) S_\kappa(s-b) S_\kappa(s-c)}}{S_\kappa (a)\, S_\kappa (b)\, S_\kappa (c)}.$$
\textbf{Proof of (B-4) :}\quad 
By (A-2), (B-1) and (B-2) we get,
\begin{eqnarray*} \sin\left({\dfrac{\alpha+\beta}{2}}\right) &=&\sin\dfrac{\alpha}{2}\, \cos\dfrac{\beta}{2} + \cos\dfrac{\alpha}{2}\, \sin\dfrac{\beta}{2}\\ 
&=&\sqrt{\dfrac{S_\kappa(s-b) S_\kappa(s-c)}{S_\kappa (b)\, S_\kappa (c)}} \sqrt{\dfrac{S_\kappa (s)\, S_\kappa(s-b)}{S_\kappa (a)\, S_\kappa (c)}}~~~~~~~~~~~~~~~~~~~~~~~~\\[1mm] &&  ~+~ \sqrt{\dfrac{S_\kappa (s)\,S_\kappa(s-a)}{S_\kappa (b)\,S_\kappa (c)}}\sqrt{\dfrac{S_\kappa(s-a)S_\kappa(s-c)}{S_\kappa (a)\,S_\kappa (c)}}~~~~~~~~~~~~~~~~~~~~~~~~~~~\\[1mm] 
&=&\sqrt{\dfrac{S_\kappa (s)\, S_\kappa(s-c)}{S_\kappa (a) \,S_\kappa (b)}} \left({\dfrac{S_\kappa(s-b) + S_\kappa(s-a)}{S_\kappa (c)}}\right)\\[1mm]
&=& \cos\left( \dfrac{\gamma}{2}\right) \, \left({\dfrac{S_\kappa(s-b) + S_\kappa(s-a)}{S_\kappa (c)}}\right)~~~~~~~~~~~~~~~~~~~~~~~~~\mbox{[by (B-2)]}
\\[1mm] 
&=& \cos\left( \dfrac{\gamma}{2}\right) \,\dfrac{2 S_\kappa\left({\dfrac{2s-a-b}{2}}\right) C_\kappa\left({\dfrac{a-b}{2}}\right)}{S_\kappa (c)}~~\mbox{[by (A-13) and (A-1)]}\\[1mm] 
&=& \cos\left( \dfrac{\gamma}{2}\right) \, \dfrac{2 S_\kappa\left( \dfrac{c}{2}\right) \, C_\kappa\left({\dfrac{a-b}{2}}\right)}{2 \,S_\kappa\left( \dfrac{c}{2}\right) \, C_\kappa\left( \dfrac{c}{2}\right) \,}~~~~~~~~~~~~~~~~~~~~~~~~~~~~~~\mbox{[by (A-4)]}\\[1mm] 
&=& \cos\left( \dfrac{\gamma}{2}\right) \, \dfrac{C_\kappa\left({\dfrac{a-b}{2}}\right)}{C_\kappa\left( \dfrac{c}{2}\right) \,} .\end{eqnarray*}
\textbf{Proof of (B-5) :}\quad Similarly, by (A-3), (B-1) and (B-2) we get,
\begin{eqnarray*} 
\sin\left({\dfrac{\alpha-\beta}{2}}\right) &=& sin\dfrac{\alpha}{2}\, \cos\dfrac{\beta}{2} - \cos\dfrac{\alpha}{2}\, \sin\dfrac{\beta}{2}
=  \cos\left( \dfrac{\gamma}{2}\right) \, \left({\dfrac{S_\kappa(s-b) - S_\kappa(s-a)}{S_\kappa (c)}}\right) \\[1mm] 
&=& \cos\left( \dfrac{\gamma}{2}\right) \,\,\dfrac{2C_\kappa\left( \dfrac{c}{2}\right) \, S_\kappa\left({\dfrac{a-b}{2}}\right)}{2 S_\kappa\left( \dfrac{c}{2}\right) \, C_\kappa\left( \dfrac{c}{2}\right) \,} \ ~~~~~~~~~~~~~\mbox{[by (A-4) and (A-14)]}\\
&=&\  \cos\left( \dfrac{\gamma}{2}\right) \,\,\dfrac{S_\kappa\left({\dfrac{a-b}{2}}\right)}{S_\kappa\left( \dfrac{c}{2}\right) \,} .~~~~~~~~~~~~~\end{eqnarray*} 
\textbf{Proof of (B-6) :}\quad 
By (A-6), (B-1) and (B-2) we get,
\begin{eqnarray*}
\cos\left({\dfrac{\alpha+\beta}{2}}\right) 
&=& \cos\dfrac{\alpha}{2}\, \cos\dfrac{\beta}{2}- \sin\dfrac{\alpha}{2}\, \sin\dfrac{\beta}{2}\\
&=& \sqrt{\dfrac{S_\kappa (s)\,S_\kappa(s-a)}{S_\kappa (b)\, S_\kappa (c)}} \sqrt{\dfrac{S_\kappa (s)\, S_\kappa(s-b)}{S_\kappa (a)\, S_\kappa (c)}}
\\[1mm]
&& ~-\, \sqrt{\dfrac{S_\kappa(s-b) S_\kappa(s-c)}{S_\kappa (b)\, S_\kappa (c)}} \sqrt{\dfrac{S_\kappa(s-a) S_\kappa(s-c)}{S_\kappa (a)\, S_\kappa (c)}}\\[1mm]
&=& \sqrt{\dfrac{S_\kappa(s-a) S_\kappa(s-b)}{S_\kappa (a) \,S_\kappa (b)}} \left({\dfrac{S_\kappa (s)\, - S_\kappa(s-c)}{S_\kappa (c)}}\right)\\
&=& \sin\dfrac{\gamma}{2}\, \left({\dfrac{S_\kappa (s)\, - S_\kappa(s-c)}{S_\kappa (c)}}\right) ~~~~~~~~~~~~~~~~~~~~~~~~~~~~~~~~\mbox{[by (B-1)]}\\[1mm]
&=& \sin\dfrac{\gamma}{2}\,\,\dfrac{2 C_\kappa\left({\dfrac{2s-c}{2}}\right) S_\kappa\left( \dfrac{c}{2}\right) \,}{2 S_\kappa\left( \dfrac{c}{2}\right) \, C_\kappa\left( \dfrac{c}{2}\right) \,}~~~~~~~~~~~~\,~~~\mbox{[by (A-4) and (A-14)]}\\[1mm]
&=&  \sin\dfrac{\gamma}{2}\,\,\dfrac{C_\kappa\left({\dfrac{a+b}{2}}\right)}{C_\kappa\left( \dfrac{c}{2}\right) \,} .
 \end{eqnarray*}
\textbf{Proof of (B-7) :} Similarly, by (A-17), (B-1), (B-2), (A-4), (A-13) and (A-1) we get 
\begin{eqnarray*}  
\cos\left({\dfrac{\alpha-\beta}{2}}\right) 
&=& \cos\dfrac{\alpha}{2}\, \cos\dfrac{\beta}{2} + \sin\dfrac{\alpha}{2}\,\sin\dfrac{\beta}{2}
= \sin\dfrac{\gamma}{2}\, \left({\dfrac{S_\kappa (s)\, + S_\kappa(s-c)}{S_\kappa (c)}}\right) ~~~~~~~~~~~~~~~~~~~~~~~~~~~~~~~~~\\[1mm]
&=& \sin\dfrac{\gamma}{2}\,\,\dfrac{2 S_\kappa\left({\dfrac{2s-c}{2}}\right)C_\kappa\left( \dfrac{c}{2}\right) \,}{2S_\kappa\left( \dfrac{c}{2}\right) \,C_\kappa
\left( \dfrac{c}{2}\right) \,}
=\sin\dfrac{\gamma}{2}\,\, \dfrac{S_\kappa\left({\dfrac{a+b}{2}}\right)}{S_\kappa\left( \dfrac{c}{2}\right) \,}.~~~~~~~~~~~~~~~~~~~~~~~~~~~~~~~~~~~~~~~~~~~~~~~~~~~~~\blacksquare  \end{eqnarray*}
\noindent \textbf{\large 6.3 Appendix C}\\[2mm]
\noindent \textbf{Proof of Proposition 2.4 :}
\begin{itemize}
\item[\textbf{(i)}] $ \boldsymbol{\Longrightarrow}$ \textbf{(ii)}\quad Let $f$ be an isometry of $M_\kappa^2$ such that $f(P)=P^\prime$, $f(Q)=Q^\prime$ and $f(R)=R^\prime$. Then, $a := d_\kappa(P,Q)=d_\kappa(f(P), f(Q))= d_\kappa(P^\prime,Q^\prime)=:a^\prime$. Similarly, $b=b^\prime $ and $c=c^\prime$.
\item[\textbf{(ii)}] $ \boldsymbol{\Longrightarrow}$ \textbf{(i)}\quad This follows immediately from Proposition 1.1.
\item[\textbf{(ii)}] $ \boldsymbol{\Longrightarrow}$ \textbf{(iii)}\quad By the Law of Cosine it follows that \\
\noindent \underline{$\kappa=0$} : 
$$\cos \alpha = \dfrac{b^2+ c^2-a^2}{2\, \,b\,c}= \dfrac{(b^\prime)^2+ (c^\prime)^2-(a^\prime)^2}{2\, \,b^\prime\,c^\prime}= \cos \alpha^\prime.$$
\noindent \underline{$\kappa \neq 0$} :$$\cos \alpha = \dfrac{C_\kappa(a)-C_\kappa(b)\, C_\kappa(c)}{\kappa\, S_\kappa(b)\,S_\kappa(c)}= \dfrac{C_\kappa(a^\prime)-C_\kappa(b^\prime)\, C_\kappa(c^\prime)}{\kappa\, S_\kappa(b^\prime)\,S_\kappa(c^\prime)}= \cos \alpha^\prime.$$

As $\alpha, \alpha^\prime \in (0,\pi)$ we get $\alpha = \alpha^\prime$.
\item[\textbf{(iii)}] $ \boldsymbol{\Longrightarrow}$ \textbf{(ii)}\quad By the Law of Cosine we have \\
\noindent \underline{$\kappa=0$} : 
$$a^2= b^2+ c^2-2\, \,b\,c\,\cos \alpha = (b^\prime)^2+ (c^\prime)^2-2\, \,b^\prime\,c^\prime\,\cos \alpha^\prime= (a^\prime)^2.$$
Since $a,^\prime >0$ we get $a= a^\prime$.\\
\noindent \underline{$\kappa \neq 0$} : 
\begin{eqnarray*}
C_\kappa(a)&=& C_\kappa(b)\, C_\kappa(c)+ \kappa\, S_\kappa(b)\,S_\kappa(c)\,\cos \alpha \\
&=& C_\kappa(b^\prime)\, C_\kappa(c^\prime)+ \kappa\, S_\kappa(b^\prime)\,S_\kappa(c^\prime)\,\cos \alpha^\prime\\
&=& C_\kappa(a^\prime)
\end{eqnarray*}
Let $$I_\kappa^\prime:= \begin{cases} (0,\pi) & \text{if $\kappa =1$,}\\
(0, \infty) & \text{if $\kappa =-1$.} \end{cases}$$ Since $a, a^\prime \in S_\kappa$ we get $a = a^\prime$.\\

Similarly we can prove that $b = b^\prime$ and $c = c^\prime$.
\item[\textbf{(iii)}] $ \boldsymbol{\Longrightarrow}$ \textbf{(iv)}\quad As (iii) $ \Longrightarrow$ (ii) we have $a=a^\prime$, $b = b^\prime$, $c = c^\prime$ and $\alpha=\alpha^\prime$. By the Law of Cosine it follows that $\beta= \beta^\prime$ and $\gamma= \gamma^\prime$.
\item[\textbf{(iv)}] $ \boldsymbol{\Longrightarrow}$ \textbf{(iii)}\quad By the Law of Sine applied to the triangles $T$ and $T^\prime$ we have,
$$\dfrac{S_\kappa(a)}{\sin \alpha^\prime}=\dfrac{S_\kappa(b^\prime)}{\sin \beta}=\dfrac{S_\kappa(c^\prime)}{\sin \gamma}.$$
and 
$$\dfrac{S_\kappa(a)}{\sin \alpha}=\dfrac{S_\kappa(b)}{\sin \beta}=\dfrac{S_\kappa(c)}{\sin \gamma}.$$
Hence,
$$\dfrac{S_\kappa(b)}{S_\kappa(c)}=\dfrac{\sin \beta}{\sin \gamma}= \dfrac{S_\kappa(b^\prime)}{S_\kappa(c^\prime)}.$$
Therefore $$S_\kappa(b)\,S_\kappa(c^\prime)= S_\kappa(b^\prime)\,S_\kappa(c).\eqno{(C-1)} $$
By the Law of Cosine we have,\\
\noindent \underline{$\kappa=0$} : 
$$\cos \gamma= \dfrac{a^2+b^2-c^2}{2\,a\,b}= \dfrac{2 \,a^2 + 2\, a\,c \, \cos \beta}{2\,a\,b}= \dfrac{a+c\, \cos \beta}{b}.$$
Similarly, 
$$\cos \gamma^\prime =
\dfrac{a^\prime+c^\prime\, \cos \beta^\prime}{b^\prime}.$$
From the hypothesis it follows that 
$$\dfrac{a+c\, \cos \beta}{b}=\dfrac{a+c^\prime\, \cos \beta}{b^\prime}.$$
Therefore,
$$a\,(b^\prime-b)= \cos \beta \, (b^\prime \,c- b\, c^\prime).$$
Thus from (C-1) it follows that $b= b^\prime$. Now, $a=a^\prime$, $b= b^\prime$ and $\gamma=\gamma^\prime$ is another form of (iii).\\
\noindent \underline{$\kappa \neq 0$} : 
\begin{eqnarray*}\cos \gamma
&=& \dfrac{C_\kappa (c)-C_\kappa(a)\, C_\kappa(b)}{\kappa\,S_\kappa(a)\, S_\kappa(b) }\\
&=&\dfrac{C_\kappa (c)-C_\kappa(a)\, \left[C_\kappa(a)\,C_\kappa(c)+ \kappa \,S_\kappa(a)\,S_\kappa(c)\, \cos \beta \right] }{\kappa\,S_\kappa(a)\, S_\kappa(b) }\\
&=&  \dfrac{\kappa\,C_\kappa (c)\,S_\kappa^2(a)+ \kappa\, S_\kappa(a)\,C_\kappa(a)\,S_\kappa(c)\, \cos \beta}{\kappa\,S_\kappa(a)\, S_\kappa(b)} ~~~~~~~~~~~~~~~~~~\mbox{[by (A-5)]}
.
\end{eqnarray*}
Since $a \in I_\kappa^\prime$, $S_\kappa(a) \neq 0$. Therefore we get
$$\cos \gamma= \dfrac{C_\kappa (c)\,S_\kappa(a)+ C_\kappa(a)\,S_\kappa(c)\, \cos \beta}{ S_\kappa(b)}.$$
Similar calculations on Triangle $T^\prime$ yields 
\begin{eqnarray*}\cos \gamma^\prime
&=& \dfrac{C_\kappa (c^\prime)\,S_\kappa(a^\prime)+ C_\kappa(a^\prime)\,S_\kappa(c^\prime)\, \cos \beta^\prime}{ S_\kappa(b^\prime)}\\
&=& \dfrac{C_\kappa (c^\prime)\,S_\kappa(a)+ C_\kappa(a)\,S_\kappa(c^\prime)\, \cos \beta}{ S_\kappa(b^\prime)}.
\end{eqnarray*}
Since $\gamma=\gamma^\prime$ we get 
\begin{eqnarray*}
& &S_\kappa(b^\prime)\left[ C_\kappa (c)\,S_\kappa(a)+ C_\kappa(a)\,S_\kappa(c)\, \cos \beta\right]\\
&=& S_\kappa(b) \left[C_\kappa (c^\prime)\,S_\kappa(a)+ C_\kappa(a)\,S_\kappa(c^\prime)\, \cos \beta \right] .
\end{eqnarray*}
That is,
\begin{eqnarray*}
& & S_\kappa(a)\,\left[ S_\kappa(b) \,C_\kappa (c^\prime)-S_\kappa(b^\prime) \,C_\kappa (c)\right] \\
&=& \kappa\, C_\kappa(a)\, \cos \beta\,\left[ S_\kappa(b) \,S_\kappa (c^\prime)-S_\kappa(b^\prime) \,S_\kappa (c)\right]. ~~~~~~~~~~~~~~~~\,~~~~~~~~~~\hfill (C-2) 
\end{eqnarray*}
As $S_\kappa(a) \neq 0$, from (C-1) and (C-2) it follows that
$$S_\kappa(b) \,C_\kappa (c^\prime)-S_\kappa(b^\prime) \,C_\kappa (c)=0.$$
That is, $$\dfrac{S_\kappa(b)}{S_\kappa(b^\prime)}= \dfrac{C_\kappa (c)}{C_\kappa (c^\prime)}.\eqno{(C-3)}$$
From (C-1) and (C-3) we get 
$$T_\kappa(c)= T_\kappa(c^\prime).$$ As $c, c^\prime \in I_\kappa^\prime$ we get $c=c^\prime$. Now, $a=a^\prime$, $c= c^\prime$ and $\beta=\beta^\prime$ is another form of (iii).\\
\item[\textbf{(ii)}] $ \boldsymbol{\Longrightarrow}$ \textbf{(v)}\quad In the proof of (ii) $ \Longrightarrow$ (iii), we showed that $\alpha = \alpha^\prime$. Similarly it can be proved that $\beta = \beta^\prime$ and $\gamma = \gamma^\prime$.\\ 
\item[\textbf{(v)}] $ \boldsymbol{\Longrightarrow}$ \textbf{(iv)}\quad \textbf{when $\boldsymbol{\kappa \neq 0}$ :}\quad Let $A, A^\prime $ denote areas of triangles $T, T^\prime$ respectively. By Proposition 2.1, 
$$A= \kappa\,\left( \pi - \alpha+\beta+\gamma\right)= \kappa\,\left( \pi - \alpha^\prime+\beta^\prime+\gamma^\prime\right)=A^\prime.$$ 
Therefore, by Proposition 2.3 it follows that
$$\dfrac{CT_\kappa\left(\frac{a}{2} \right)\,CT_\kappa\left(\frac{b}{2} \right)+ \kappa\, \cos \gamma }{\sin \gamma}=\dfrac{CT_\kappa\left(\frac{a^\prime}{2} \right)\,CT_\kappa\left(\frac{b^\prime}{2} \right)+ \kappa\, \cos \gamma^\prime }{\sin \gamma^\prime}.$$
Since $\gamma= \gamma^\prime$ we get 
$$CT_\kappa\left(\frac{a}{2} \right)\,CT_\kappa\left(\frac{b}{2}\right) =CT_\kappa\left(\frac{a^\prime}{2} \right)\,CT_\kappa\left(\frac{b^\prime}{2}\right) .\eqno{(C-4)}$$
Similarly we have
$$CT_\kappa\left(\frac{a}{2} \right)\,CT_\kappa\left(\frac{c}{2}\right) =CT_\kappa\left(\frac{a^\prime}{2} \right)\,CT_\kappa\left(\frac{c^\prime}{2}\right) ,\eqno{(C-5)}$$
and
$$CT_\kappa\left(\frac{b}{2} \right)\,CT_\kappa\left(\frac{c}{2}\right) =CT_\kappa\left(\frac{b^\prime}{2} \right)\,CT_\kappa\left(\frac{c^\prime}{2}\right) .\eqno{(C-6)}$$
Multiplying (C-4) and (C-5) we get
$$CT_\kappa^2\left(\frac{a}{2} \right)\,CT_\kappa\left(\frac{b}{2}\right) \,CT_\kappa\left(\frac{c}{2}\right) =CT_\kappa^2\left(\frac{a^\prime}{2} \right)\,CT_\kappa\left(\frac{b^\prime}{2}\right) \,CT_\kappa\left(\frac{c^\prime}{2}\right) .$$
As $CT_\kappa(x) \neq 0$ for $x \in I_\kappa^\prime$, by (C-6) we get 
$$CT_\kappa^2\left(\frac{a}{2}\right) =CT_\kappa^2\left(\frac{a^\prime}{2}\right) .$$
That is, $$T_\kappa^2\left(\frac{a}{2}\right) =T_\kappa^2\left(\frac{a^\prime}{2}\right) .\eqno{(C-7)}$$
As $a, a^\prime \in I_\kappa^\prime$, $\frac{a}{2}, \frac{a^\prime}{2} \in I_\kappa^{\prime\prime}$ where $I_\kappa^{\prime\prime}=\begin{cases} (0, \frac{\pi}{2})& \text{if $\kappa=1$}\\
(0, \infty) &\text{if $\kappa =-1$}.                                              \end{cases}$
Therefore, $$T_\kappa (a), ~T_\kappa(a^\prime )>0$$ and hence from (C-7) we get $T_\kappa (a)= T_\kappa(a^\prime )$. Finally, since $T_\kappa$ is an injective function on $I_\kappa^{\prime\prime}$ we get $a= a^\prime$. \\
Similarly it can be shown that $b= b^\prime$ and $c= c^\prime$. Thus, in fact we have (v) $\Longrightarrow$ (ii).
\end{itemize} This completes the proof of Proposition 2.4.\hfill $\blacksquare$ \\[3mm]

\noindent {\sf A. R. Aithal}\\
\begin{small}
Retired, Department of Mathematics\\
University of Mumbai, Mumbai-400\,098.\\
aithal86@gmail.com\\[3mm]
\end{small}
{\sf Anisa Chorwadwala} \\
\begin{small}
Indian Institute of Science Education and Research Pune, \\
Dr Homi Bhabha Road, Pashan, Pune-411008. \\
anisa@iiserpune.ac.in  
\end{small}
\end{document}